\documentclass[12pt,
 oneside
]{amsart}
\usepackage{amscd,amsmath, wasysym}
\title[]{Common boundary values of holomorphic functions for
two-sided complex structures}
\author[]{Florian Bertrand \and  Xianghong Gong$^*$ \and Jean-Pierre Rosay}

 \address{Department of Mathematics,
 University of Wisconsin, Madison, WI 53706, U.S.A.}
\email{bertrand@math.wisc.edu, gong@math.wisc.edu, jrosay@math.wisc.edu.}

  \thanks{$^*$ Research  supported in part by NSF grant DMS-0705426.}

\newcommand{\dist}{\operatorname{dist}}

\newtheorem{thm}{Theorem}[section]
\newtheorem{cor}[thm]{Corollary}
\newtheorem{prop}[thm]{Proposition}
\newtheorem{lemma}[thm]{Lemma}

\theoremstyle{definition}

\newtheorem{exmp}[thm]{Example}

\newtheorem{rem}[thm]{Remark}

\renewcommand{\th}[1]{\begin{thm}\label{#1}}
\newcommand{\eth}{\end{thm}}
\newcommand{\co}[1]{\begin{cor}\label{#1}}
\newcommand{\eco}{\end{cor}}
\renewcommand{\le}[1]{\begin{lemma}\label{#1}}
\newcommand{\ele}{\end{lemma}}
\newcommand{\pr}[1]{\begin{prop}\label{#1}}
\newcommand{\epr}{\end{prop}}

\newcommand{\ga}{\begin{gather}}
\newcommand{\ega}{\end{gather}}
\newcommand{\gan}{\begin{gather*}}
\newcommand{\egan}{\end{gather*}}
\newcommand{\al}{\begin{align}}
\newcommand{\eal}{\end{align}}
\newcommand{\aln}{\begin{align*}}
\newcommand{\ealn}{\end{align*}}
\newcommand{\eq}[1]{\begin{equation}\label{#1}}
\newcommand{\eeq}{\end{equation}}

\newcommand{\pdoz}{\partial_{\overline z}}

\newcommand{\ci}{~\cite}

\newcommand{\f}[2]{\frac{#1}{#2}}

\newcommand{\df}{\overset{\text{def}}{=\!\!=}}
\newcommand{\D}{\mathbb{D}}
\newcommand{\cc}{{\bf C}}

\newcommand{\rr}{{\bf R}}

\newcommand{\rtn}{{\bf R}^{2n}}

\newcommand{\ov}{\overline}
\newcommand{\ord}{\operatorname{ord}}

\newcommand{\RE}{\operatorname{Re}}
\newcommand{\IM}{\operatorname{Im}}
\newcommand{\dbar}{\overline\partial}

\newcommand{\cL}{\mathcal}

\newcommand{\I}{\operatorname{I}}
\newcommand{\oS}{ S}
\newcommand{\T}{ T}
\newcommand{\Ct}{\mathcal C}

\newcommand{\B}{\mathcal C}
\newcommand{\hB}{\hat{\mathcal C}}

\newcommand{\Ht}{\mathcal H}

\newcommand{\all}{\alpha}

\newcommand{\gaa}{\gamma}

\newcommand{\del}{\delta}
\newcommand{\Del}{\Delta}
\newcommand{\var}{\varphi}
\newcommand{\e}{\epsilon}
\newcommand{\om}{\omega}
\newcommand{\Om}{\Omega}

\newcommand{\thh}{\theta}
\newcommand{\la}{\lambda}
\newcommand{\ta}{\tau}
\newcommand{\fta}{\varsigma}

\newcommand{\nm}[2]{|#1|_{#2}}

\newcommand{\pd}{\partial}

\newcommand{\yt}{\frac{1}{2}}

\newcommand{\re}[1]{(\ref{#1})}
\newcommand{\rea}[1]{$(\ref{#1})$}
\newcommand{\rl}[1]{Lemma~\ref{#1}}

\newcommand{\rp}[1]{Proposition~\ref{#1}}
\newcommand{\rt}[1]{Theorem~\ref{#1}}

\newcommand{\supp}{\operatorname{supp}}

\newcommand{\db}{\dbar}

\newcounter{pp}
\newcommand{\bpp}{\begin{list}{$\hspace{-1em}\alph{pp})$}{\usecounter{pp}}}
\newcommand{\epp}{\end{list}}

\newcounter{ppp}
\newcommand{\bppp}{\begin{list}{$\hspace{-1em}(\roman{ppp})$}{\usecounter{ppp}}}
\newcommand{\eppp}{\end{list}}

\begin{document}
\begin{abstract}
Let $\Omega_1,\Omega_2$ be two disjoint open  sets in $\mathbf
C^n$ whose boundaries share a smooth real hypersurface $M$ as
relatively open subsets. Assume that $\Omega_i$ is equipped with a
complex structure $J^i$ which is    smooth up to $M$. Assume that
the operator norm $\|J^2-J^1\|<2$ on $M$.
 Let $f$ be a continuous  function  on
the union of $\Omega_1,\Omega_2, M$. If $f$   is holomorphic
with respect to both  structures in the open sets, then $f$ must be smooth on the union of $\Omega_1$
with $M$.
Although the result as stated is far more
 meaningful for integrable structures, our methods make it much more natural to deal with the  general
  almost complex structures without the integrability condition.
 The result is therefore proved
  in the framework of almost complex structures.
%
\end{abstract}

 \maketitle

\tableofcontents

\setcounter{section}{0}
\setcounter{thm}{0}\setcounter{equation}{0}
\section{Introduction}\label{sec1}

The title of the paper is suggested by the following result.
\pr{combv}    Let $0<\all<1$ and $k\geq0$ be an integer.
Let $\gaa$ be an embedded curve in $\cc$ of class
$\cL C^{k+1+\all}$. Let $\Om_1, \Om_2$ be disjoint open subsets of $\cc$. Suppose that
both boundaries $\pd\Om_1$,   $\pd\Om_2$ contain
$\gaa$ as  relatively open subsets. Assume that
 $a_i\in\cL C^{k+\all}(\Om_i\cup\gaa)$   satisfy $
|a_i(z)|<1$ on $\Om_i\cup\gaa$.
Let $f$ be a continuous function on $ \Om_1\cup\gaa\cup \Om_2$   satisfying
$$
\pdoz f+a_i \pd_zf=0\quad \text{on $\Om_i$},\quad i=1,2.
$$
Then   $f$ is in  $\cL C_{loc}^{k+1+\all}
(\Om_1\cup\gaa)\cap\cL C_{loc}^{k+1+\all}
( \Om_2\cup\gaa)$.
\epr
Here $  \cL C_{loc}^{k+\all}
( \Om\cup\gaa)$ denotes the set of functions $f\in \cL C^{k}
( \Om\cup\gaa)$ whose $k$-th order derivatives are in $\cL C^\all(K)$ for each compact subset of $\Om\cup\gaa$.
The result does not hold in general for harmonic functions
due to a jump formula
 for single-layer distributions. A more general result for non-homogeneous equations
  is in \rp{combv+}.
Our next result is in higher dimension.
%
%
%
%
%
%
\th{combvn}
 Let   $k\geq  4,   n\geq2$ be  integers and let $0<\all<1$.  Let $\Om_1,\Om_2$ be disjoint open subsets of $\cc^n$
 such that   both boundaries $\pd\Om_1$, $\pd \Om_2$ contain a smooth real hypersurface
$M$ of class $\cL C^{k+1+\all}$ as  relatively open subsets.   For $i=1,2$, let $J^i$ be an almost  complex
structure of class $\cL C^{k+\all}$ on $\Om_i\cup M$. Suppose that on $M$
 the operator norm $\|J^2-J^1\|<2$.
Let $f$ be a continuous function on $\Om_1\cup M\cup \Om_2$. Suppose that  for $i=1,2$, $1\leq j\leq n$,
  $ (\pd_{x_j}+\sqrt{-1}J^i\pd_{x_j})f$ and $ (\pd_{y_j}+\sqrt{-1}J^i\pd_{y_j})f$, defined on $\Om_i$,  extend to
 functions in $ \cL C^{  k }( \Om_i\cup M)$.
 Then   $f$ is of class $\cL C^{k-3+\beta}(\Om_1\cup M)$ for all $\beta<1$.
In particular, $f\in\cL C^\infty(\Om_1\cup M)$ when $k=\infty$.
\eth
Notice that no assumption is made on convexity of $M$ with respect to either of the  almost complex structures.
The definition of almost complex structures is in section~\ref{sec3} and
a general result is in section~\ref{sec4}.

We would like to mention that
the interior regularity of $f$ for integrable almost complex structures
 is ensured by the well-known Newlander-Nirenberg theorem~\ci{NNfise} (see also
Nijenhuis-Woolf~\ci{NWsith} and Webster~\ci{Weeini}).
There are   results on   Newlander-Nirenberg theorem for   pseudoconvex
domains with boundary by Catlin~\ci{Caeiei} and Hanges-Jacobowitz~\ci{HJeini}.
See   earlier work of
 Hill~\ci{Hieiei} on failure  of Newlander-Nirenberg type theorem
with boundary.

\smallskip

We now observe how  the
 common boundary values 
arise  in  the Cauchy-Green
  operator for $\dbar$ in $\cc$.
Let $\Om$ be a bounded domain in $\cc$ with $\pd\Om\in\cL C^{k+1+\all}$.
Seeking coordinates $z+f(z)$ to transform  $\pdoz+a\pd_z$
into a multiple of $\pdoz$ leads to the   equation
\eq{diffeq}
\pdoz f+a(z)\pd_zf+b(z)=0, \quad z\in\Om,
\eeq
where $a\in\cL C^{k+\all}(\ov\Om)$, and $b$ is either $a$ or   a function of the same
kind.
 To solve it, one   considers the integro-differential equation
\eq{tdiffeq}
f(z)+T(a\pd_zf)(z)+Tb(z)=0, \quad z\in\Om.
\eeq
Here $T=T_\Om$ is the Cauchy-Green operator
$$
Tf(z)=\f{1}{\pi}\int_{\Om}\f{f(\zeta)}{z-\zeta }\, d\xi\,
d\eta.
$$
 The
equation \re{tdiffeq} is equivalent
to \re{diffeq} and an extra equation
\eq{0=cf}
 \int_{\pd\Om}\f{f(\zeta)}{\zeta-z}\, d\zeta=0, \quad z\in\Om.
\eeq
When $f\in\cL C^0(\pd\Om)$,   the jump formula implies that
\re{0=cf} is equivalent to $f$ being the boundary value of
a  function which is holomorphic on $\Om'=\cc\setminus \ov\Om$,
continuous on $\ov{\Om'}$, and vanishing at $\infty$.
(See~\rl{teqs1}  for details.)

As an application of \rp{combv}, we will prove the following.
\th{regid-}
Let $0<\all<1$ and let
  $\Om\subset\cc$ be a   bounded domain  with
$\cL C^{1+\all}$ boundary.
Let $a,b\in\cL C^{\all}(\ov\Om)$. There exists $\e_\all>0$ such that
if
   $\|{a}\|_{\all}<\e_\all$, then     \rea{tdiffeq}  admits a unique
   solution $f\in\cL C^{1+\all}(\ov\Om)$.
   Assume further that  $a,b\in\cL C^{k+\all}(\ov\Om)$ and
    $\pd\Om\in \cL C^{k+1+\all}$ for an integer $k\geq0$.
 Then $f\in\cL C^{k+1+\all}(\ov\Om)$. Consequently,
the linear map $\I+Ta\pd_z$ from $\cL C^{k+1+\all}(\ov\Om)$ into itself
has a bounded inverse.
\end{thm}

The structures $\pdoz$ and $\pd_z$ show that \rt{combvn} fails for $\|J^2-J^1\|=2$ and $M\colon \IM z_1=0$.
We   expect that the regularity of  $f$ in \rt{combvn} can be improved.
The loss of derivatives is due to an essential use of the
Fourier transform.
For this reason, we will present two alternative proofs for  the one-dimensional
case, with one producing the sharp  result.


\medskip

We want to mention two open problems in addition to the regularity issue
 mentioned above.  The first problem is on the vector-valued version of \re{tdiffeq}.
 The second is  concerned with  non-linear integro-differential equations
 arising from differential equations of $J$-holomorphic curves; in fact,
   the integrable case remains to be studied.

\medskip
\noindent{\bf Problem A.} {\it Let $m\geq2, 0<\all<1$ and let $D$ be
a bounded domain in $\cc$ with $\cL C^\infty$ boundary. Let
$A=(a_{jk})\in\cL C^\infty(\ov D)$ be an $m\times m$ matrix  with
sufficiently small $\cL C^\all$ norm   on $\ov D$.  Does $ \I+T_D  A
{\pd_z }  \colon [ \cL C^{k+\all}(\ov D)]^m\to [\cL C^{k+\all}(\ov
D)]^m$ have a  bounded inverse for
 all positive integer $k$?}

\smallskip

\noindent{\bf Problem B.}  {\it  Let $D$ be a bounded domain in $\cc$ with $\cL C^\infty$ boundary.
 Let $\Om$ be a domain in $\cc^n$ with   $n\geq1$.
Let $A$ be an $n\times n$ matrix of $ \cL C^{\infty}$ functions on $\Om$.
  Suppose that the operator norm $\|A(z)\|$  is less than $1$
  on $\Om$.
Suppose that a continuous map $u\colon\ov D\to \Om$ satisfies
\eq{utdv}
u+T_D (A(u)\ov{\pd_zu})=T_Dv
\eeq
on $D$. Here  $v$ is a $\cL C^\infty$ map from $ \ov D$ into $ \cc^n$.
Is $u\in\cL C^\infty(\ov D)$?}

  Note that  the interior
 regularity of $u$ is in work of  Nijenhuis-Woolf~\ci{NWsith}.  When the $\cL C^{1+\all}$
 norm of $A$ is sufficiently small, the existence and uniqueness of solutions $u$ to \re{utdv}
 is ensured.
 Problems A and B can be reformulated in terms of two differential equations
 on $D$ and its complement. Indeed, by   \rl{teqs1}   a continuous map $u\colon\ov D\to\Om$
  satisfies \re{utdv} if and only if it extends to a continuous map $u$ from $\cc$ into itself that vanishes
   at $\infty$ and satisfies
 $\pdoz u+A(u)\ov{\pd_zu}=v$ on $D$ and $\pdoz u=0$ on $\cc\setminus\ov D$.

\setcounter{thm}{0}\setcounter{equation}{0}
\section{Inverting $\I+TA\ov{\pd_z}$
}\label{sec2}

In this section, we will recall     estimates on the Cauchy-Green operator $T$
and $ \pd_zT$. We will discuss the inversion
of $\I+TA\pd_z$ in spaces of higher order derivatives when $A$
has a   small    $\cL C^\all$ norm.
When $A$ has compact support, it is easy to bound     inverses of  $\I+TA\pd_z,
\I+TA\ov{\pd_z}$.
We will show in  section~\ref{sec5} that $\I+TA\pd_z$ is indeed invertible when $A$ is a suitable
scalar  function.

 Throughout the paper, when a parameter set $P$ is involved in $\Om\times P$,
   $\Om$ is a bounded open set in a euclidean space and
$P$ is the closure of a  bounded open set in a euclidean space. We assume that
      two points $a,b$ in $\ov\Om\times  P$ can be connected
by a smooth curve in $\ov\Om\times   P$ of length at most $C|b-a|$.

We will need
  spaces of functions with parameter.
  The usual norm on $\cL C^{k+\all}(\ov\Om\times   P)$ is denoted by   $|\cdot|_{k+\all}$.
Following  \ci{NWsith}, for     integers $k,j\geq0$ we define $\hB^{k+\all,j}(\ov\Om,  P)$
to be the set of
functions $f(z,t)$ such that for all $i\leq j$, $\pd_t^if\in \cL C^{k}(\ov\Om\times P)$
and
$$
\|f\|_{k+\all,j}=\max_{0\leq i\leq j}|\pd_t^if(\cdot,t)|_{k+\all}<\infty.
$$
  Define $\hB^{\infty,j}(\ov\Om,  P)=\bigcap_{k=1}^\infty\hB^{k,j}(\ov\Om,  P)$.
Throughout the paper, $k$ is a nonnegative integer, and $0<\all<1$.
To simplify   notation, the parameter set
$P$ will not be indicated sometimes.

Let $\Om$ be a bounded domain in $\cc$. The $\db$ solution operator $T$ and $S=\pd_zT$ are
\ga
 \label{tfzt}
\T f(z) =\f{1}{\pi}\int_{\Om}\f{f(\zeta)}{z-\zeta}\,
d\xi\, d\eta,\quad \oS f(z)  =-\f{1}{\pi} \, p.v.\int_{ \Om}
\f{f(\zeta)}{(z-\zeta)^2 }\, d\xi\,
d\eta.
\end{gather}
It is well-known that $\pdoz \T$ is the identity on $L^p(\Om)$ when
$p>2$.
When $f\in\cL C^\all(\ov D)$ and $\pd\Om\in\cL C^{1+\all}$, one has
\ga
\label{sfzt} \oS f(z) 
 =-\f{1}{\pi} \int_{\Om}\f{f(\zeta)-f(z)}{(z-\zeta)^2 }\, d\xi\,
d\eta-\f{f(z)}{2\pi i}\int_{\pd\Om}\f{d\ov\zeta }{\zeta-z}.
\end{gather}
 If $f$ has compact support in $\Om$,
or if $f\in\cL C^{k+\alpha}(\ov{\Om})$  and $\pd\Om\in\cL C^{k+1+\all}$,
then $\T$, $\oS $ satisfy
\eq{beve}
 |Tf |_{k+1+\alpha }\leq C_{k+1+\alpha} |f |_{k+\all }, \quad
 |\oS f |_{k+\alpha}\leq C_{k+1+\alpha} |f |_{k+\all }.
\eeq
See Bers~\ci{Befise} and Vekua~\ci{Vesitw} (p.~56).  The above estimates for domains with
 parameter will be derived in
section~\ref{sec3+}. It is known that
\eq{dzsf}
\pd_zSf=S\pd_zf, 
\quad \pdoz Sf=\pd_zf,
\eeq
 where the first identity needs $f$ to have compact support in $\Om$.

For $f\in\hB^{k+\all,j}(\ov\Om, P)$,    define $Tf(z,t), Sf(z,t)$ by \re{tfzt}-\re{sfzt}
by  fixing $t$. %
\le{tspa} Let $\Om\subset\cc$ be a bounded domain with  $
\pd\Om\in\cL C^{k+1+\all}$. Then
\ga
\nonumber 
T\colon\hB^{k+\all,j}(\ov\Om,P)\to\hB^{k+1+\all,j} (\ov\Om,P), \quad
\|Tf\|_{k+1+\all,j}\leq C_{k+1+\all}\|f\|_{k+\all,j},\\
S\colon\hB^{k+\all,j}(\ov\Om,P)\to\hB^{k+\all,j} (\ov\Om,P), \quad
\|Sf\|_{k+\all,j}\leq C_{k+1+\all}\|f\|_{k+\all,j}.
\label{beve++}
\end{gather}
\ele
\begin{proof}  By \re{sfzt}, we get $S(\hB^{k+\all,j})\subset
\hB^{k,j}$.
We can verify  that $\pd_tS=S\pd_t$ on $\hB^{\all,j}$ for $j\geq1$. Thus $S(\hB^{k+\all,j})
\subset \hB^{k+\all,j}$ by \re{beve}.

The Cauchy kernel is integrable. So $T(\hB^{0,j}(\ov\Om\times P))
\subset\cL \hB^{0,j}(\ov\Om\times P)$. Also $\pd_tT=T\pd_t$ on $\hB^{0,j}$ for $j\geq1$.
The rest of assertions follows from $\pd_{z}T=S$ and $\pdoz T=\I$.
\end{proof}

By an abuse of notation, we define
$
\ov{\pd_z}f=\ov{\pd_zf}.
$
\le{isa} 
Let $\Om$ be a bounded domain in $\cc$.
Let $A\in  \hB^{k+\all,j}(\ov\Om, P )$ be an $m\times m$ matrix.
There exists     $\e_\all$ which   depends only on
$\all$ and satisfies the following.
\bppp
\item If $\pd\Om\in\cL C^{1+\all}$ and
 $\nm{A}{\all,0} < \e_{\all}$, then
$$
\I+TA\pd_z,\  \I+TA\ov{\pd_z}\colon \bigl[\hB^{1+\all,j}(\ov\Om )\bigr]^m \to
\bigl[\hB^{1+\all,j}(\ov\Om )\bigr]^m
$$
 have bounded inverses.
\item
If $A(\cdot,t)$ have compact support
in $\Om$ for all $t\in P$ and $\nm{A}{\all,0} < \e_{\all}$, then
$$
\I+TA\pd_z,\  \I+TA\ov{\pd_z}\colon \bigl[\hB^{k+1+\all,j}(\ov\Om )\bigr]^m \to
\bigl[\hB^{k+1+\all,j}(\ov\Om )\bigr]^m
$$
 have bounded inverse.
  \eppp
\end{lemma}
\begin{proof}  To be made precise, when $S$ operates on functions with compact support,
it commutes with  $\pd_t, \pd_z, \pdoz$  somewhat.   However,
 differentiating the operator product $(SA)^n$ requires
 counting
terms  efficiently as  $n$ tends to $\infty$.

(i).
Fix $0<\thh<1/2$. Note that
\gan 
\|fg\|_{k+\all,j}\leq C_{k,j}\|f\|_{k+\all,j}\|g\|_{k+\all,j}.
\end{gather*}
By \re{beve++}, we have $\|SA\|_{\all,0}\leq  C_\all'\|A\|_{\all,0}$.
Thus,
$$
\|(SA)^n\|_{\all,0}\leq  (C_\all\|A\|_{\all,0})^{n}
 \leq  \thh^{n}
$$
if $\nm{A}{\all,0}$ is sufficiently small.
Note that
$$
(TA\pd_z)^{n}=TA(SA)^{n-1}\pd_z.
$$
Let $L=\I+\sum_{n=1}^\infty (-1)^{n}TA(SA)^{n-1}\pd_z$.
Then
$$
\|TA(SA)^{n-1}\pd_zf\|_{1+\all,0}\leq C_\all\theta^{n-1}\|A\|_{\all,0}\|f\|_{1+\all,0}.
$$
This shows that
for  $f\in \hB^{1+\all,0}$,
$
\sum_{n=0}^\infty  (-1)^{n}TA(SA)^{n-1}\pd_zf
$
converges to $Lf\in\hB^{1+\all,0}$. Moreover,
$
\|Lf\|_{1+\all,0}\leq C\|f\|_{1+\all,0}.
$
It is straightforward that $L(\I+TA\pd_z)$ and $ (\I+TA\pd_z)L$
 are the identity on $\hB^{1+\all,0}$. This verifies (i) for $j=0$. The case of $j>0$
 will follow from the argument in (ii) below, by using $\pd_tT=T\pd_t$, $\pd_tS=S\pd_t$.

(ii). We need to show that $\sum\|( \oS A)^n\|_{k+\all,j}$ converges
when $A$ has compact support in $\Om$.
Denote by $C_{k+\all,j}$ a constant depending only on $k,j$,
and $\|A\|_{k+\all,j}$.
By \re{dzsf} and $\pd_tS=S\pd_t$, we
can write
$$
\pd SA=\tilde S\tilde \pd A.
$$
Here $\tilde S$ is  either $S$ or $\I$, and $\pd,\tilde \pd$
are   of form $\pd_z,\pdoz,\pd_t$.
Denote by $\pd^K$   derivatives in $z$, $\ov z$. Then
$\pd (SA)^n$ equals a sum of terms of the form
$$
S_{m_1}(\pd^{K_1} A)\cdots  S_{m_n}(\pd^{K_{n}}A)\pd^{K_{n+1}}, \quad |K_1|+\cdots+|K_{n+1}|=1.
$$
Here $S_{m_i}$ is either
$S$ or $\I$; in particular, $\|S_{m_i}\|_{k+\all,j}\leq C_{k+\all,j}$ for
all $m_i$.
The sum has at most $n+1$ terms. Thus $\pd^K\pd_t^J(SA)^{n}$ is a sum of
at most $(n+1)^{|K|+|J|}$
terms of
\eq{sKJ}
 S_{m_1}(\pd^{K_1}\pd_t^{J_1}A)\cdots S_{m_n}(\pd^{K_n}\pd_t^{J_n}A)S_{m_n}\pd^{K_{n+1}}\pd_t^{J_{n+1}}.
\eeq
Assume that $|K|\leq k, |J|\leq j$ and $n>k+j$. With $C_\all\geq1$,
\al\nonumber 
|\pd^K\pd_t^J((SA)^nf)|_{\all,0}
&\leq (n+1)^{k+j}
C_\all^n (1+|A|_{k+\all,j})^{k+j}|A|_{\all,0}^{n-k-j} \|f\|_{k+\all,j}\\
&\leq
(n+1)^{k+j}C_{k+\all,j}\|f\|_{k+\all,j}\theta^{n-k-j}.
\nonumber
\end{align}
This shows that $\|(SA)^n\|_{k+\all,j}\leq C_{k+\all,j}' (n+1)^{k+j}\theta^{n-k-j}$. Hence
$$\|TA(SA)^{n}\pd_z\|_{k+1+\all,j}\leq
C_{k+\all,j} (n+1)^{k+j}\theta^{n-k-j}.$$
We conclude that $\|(\I+TA\pd_z)^{-1}\|_{k+1+\all,j}<\infty$.

The proof for $\I+TA\ov{\pd_z}$  is obtained by minor changes. Indeed, with
 $Cf=\ov f$,  we write $
(TA\ov{\pd_z})^n=TAC(SAC)^{n-1}\pd_z.$
Now, $\pd_zC=C\pdoz$ and $\pdoz C=C\pd_z$.  We may assume that $t$ are
real variables. So $\pd_tC=C\pd_{t}$. Thus $\pd^K\pd_t^J(SAC)^n$ is a sum of at most $(n+1)^{|K|+|J|}$
terms of
$$
S_{m_1}(\pd^{K_1}\pd_t^{J_1}A)
C\cdots S_{m_n}(\pd^{K_n}\pd_t^{J_n}A)C\pd^{K_{n+1}}\pd_t^{J_{n+1}}.
$$
Substitute the above for \re{sKJ}. The remaining argument follows easily.
\end{proof}

We need a simple version of Whitney's extension theorem with parameter.
\le{whit} Let $N$ be a positive integer or $ \infty$, $0\leq j<\infty$,
 and   $0\leq\all<1$. Let $e_k$ be a sequence of positive numbers.
  Let $f_I\in\hB^{N -1-|I|+\all,j}(\rr^n,  P)$ for $0\leq |I|<N$ with $I=(i_1,\ldots, i_m)$.
Assume  that all $f_I(\cdot,t)$ have
  support contained in a compact subset $K$  of the unit ball $B$.
  There exists
  $Ef\in\hB^{ N-1+\all,j}(\rr^{n}\times\rr^m,  P)$
  such that $\pd^I_yEf(x,0,t)=f_I(x,t)$. 
   Moreover, $Ef(\cdot, t)$
  have compact support in the unit ball
  of $\rr^{n}\times\rr^{m}$ and 
\eq{exm}
\|Ef\|_{k+\all,j}\leq \e_k+C_{N,k, K}\sum_{|I|\leq k}\|f_I\|_{k-|I|+\all,j}, \quad 0\leq k<N.
\eeq
Here $E,  C_{N,k,K}$ are independent of $j$.
\ele
\begin{proof} When $N=\infty$, take   $g\in\cL C^\infty_0(\Del_\del^m)$
with $g(y)-1$ vanishing to infinity order at $0$. Set
 \eq{ymge-}
Ef(x,y,t)= \sum \f{y^I}{I!}f_I(x,t) g( \del_{|I|}^{-1} y).
 \eeq
For the above to converge, choose $\del_i$ which  decrease  to $0$
so rapidly that for all $t$
\gan 
\|b_I\|_{i-1,j}\leq \del_{i}^{1/2}, \quad
b_{I}(x,y,t)=y^Ig(\del_{i}^{-1} y) f_I(x,t), \quad   i=|I|\geq1.
\end{gather*}
Thus,  $Ef\in\hB^{k,j}$ for all $k$. Note that we can choose
$\del_k$ so small that $\|Ef\|_{k,j}\leq e_k+C_{k}\sum_{|I|\leq k}\|f_I\|_{k,j}$ for
all $k$.

Let  $N<\infty$.
Let $\phi$ be a smooth function on $\rr^n$ with support   in $B_\del^n$ for a small $\del$. Also $\int_\rr\phi(y)\, dy=1$.
We extend $f$ one dimension at a time. Assume that $m=1$.
  We need to modify   extension \re{ymge-}.  We   replace $y^if_i(x,t)$ by $y^ig_i(x,y,t)$
to achieve the   $\hB^{N-1+\all,j}$ smoothness. We also need the  correct $i$-th
  $y$-derivative  of $y^ig_i(x,y,t)$  due to the presence of  $y^lg_l(x,y,t)$ for
  $l<i$.  With $a_i\in\hB^{N-1-i,j}$ to be determined,   consider
$$
g_i(x,y,t)=\int_{\rr^n} a_i(x-yz,t)\phi(z)\, dz.
$$
  Fix $t$. We first show that
  $y^ig_i(x,y,t)$ is of class $\cL C^{N-1+\all,j}$. Since it is $\cL C^\infty$
for $y\neq0$, it suffices to extend  its
  partial derivative  of  order $<N$ on $y\neq0$    continuously to
   $\rr^n\times\rr$, as the extensions are clearly independent of
  the order of differentiation. By the product rule, this amounts to
extending $y^{i-l}\pd^Ig_i$ for $|I|\leq N-1-l$.
When $|I|\leq N-1-i$, we take $I_1=I$ and $I_2=0$. Otherwise,  $I=I_1+I_2$ with $|I_1|=N-1-i$.
Then
$$
\pd^{I_1}  g_i(x,y,t)=\sum_{|L|=|I_1|}\int \pd^La_i(x-yz,t)\phi_{I_1L}(z)\, dz
$$
for some $\phi_{I_1L}$ with   support in $B_\del^n$.
When $y\neq0$, change   variables and take    derivative $\pd^{I_2}$. We get
$$
\pd^{I}  g_i(x,y,t)=\sum_{|L|=|I_1|}
\int\f{1}{y^{ |I_2|+n}}\pd^La_i(z)
\tilde\phi_{I_1I_2L}\left(\f{x-z}{y}\right)\, dz.
$$
Change  variables again. We get
$$
y^{ |I_2|}\pd^{I}g_i(x,y,t)=\sum_{|L|=|I_1|}\int\pd^La_i(x-yz,t)\tilde
\phi_{I_1I_2L}(z)\, dz.
$$
The right-hand side and its derivatives in $t$ of order at most $j$ are
 clearly continuous functions. Since $|I_2|\leq (N-1-l)-(N-1-i)=i-l$,
then $y^{i-l}\pd^Ig_i$ extends
continuously to $\rr^n\times\rr$.
Take derivative in parameter $t$ and  compute the H\"older ratio in $x,y$. We
get
$$
\|b_{i,lI}\|_{\all,j}\leq C_{N,K}\|a_i\|_{N-1-i+\all,j},
\quad b_{i,lI}(x,y,t)=y^{i-l}\pd^{I}g_i(x,y,t).
$$
By the product rule,   at $y=0$
$$
\pd_y^{i}(y^{ i}g_i(x,y,t))=i!a_i(x,t),\quad   \pd_y^{l}(y^{ i}g_i(x,y,t))=0,\quad l<i.
$$
Starting with $a_0=f_0$, inductively we find $a_i\in\hB^{N-1-i+\all,j}$ such that
 $$
Ef(x,y,t)= \sum_{ i<N} \f{y^i}{i!}g_i(x,t) g( \del^{-1} y)
$$
  satisfies $\pd_y^iEf(x,y,t)=f_i(x,t)$
for $i=0$, $1, \ldots$, $N-1$.  The estimate  \re{exm} is immediate
when $\del$ is sufficiently small.

For $m>1$, suppose that we have found extensions $\tilde f_{i}\in\hB^{N-1-i+\all,j}(\rr^n\times\rr^{m-1},P)$
such that $\pd_{y'}^{I'}\tilde f_{i}=f_{I'i}$ at $y'=0$ for all $|I'|<N-i$ and
\eq{exm1}
 \|\tilde f_{i}\|_{k-i+\all,j}\leq e_{N, K}'+C_{N,K}\sum_{|I'|\leq k-i}\|f_{I'i}\|_{k-|I'|-i+\all,j},\quad k<N
\eeq
with $e_{N, K}'>0$ to be determined.
Assume further that $\tilde f_{i}(\cdot, t)$ have  support in a compact subset $K'$ of the unit ball of $\rr^{n+m-1}$,
where $K'$ depends only on $K$.
Using the one-dimensional   result again, we get $Ef\in \hB^{N-1+\all,j}(\rr^n\times\rr^{m},P)$ with
 compact support in the unit ball of $\rr^{n+m}$. Furthermore,   $\pd_{y_n}^{i}Ef=\tilde f_{i}$  at $y_n=0$
 and
\aln
\|Ef\|_{k+\all,j}\leq e_{N,K}'+C_{N,K}'\sum_{0\leq i\leq k}\|\tilde f_{i}\|_{k-i+\all,j}, \quad k<N.
\end{align*}
  Let $e_{N,K}'$   be sufficiently small.
Combining with \re{exm1} yields \re{exm}.
\end{proof}

  The above proof
for non-parameter case   is in~\ci{Honize} (pp.~16 and 18).
When $f$ is   defined on  $y_n\leq 0$ with $\pd_{x_n}^kf=f_k$ on $y_n=0$, the above extension
$Ef$ can be replaced by $f$ on $y_n\leq0$. The same conclusions on $Ef$ hold.
   Seeley~\ci{Sesifo}  has a linear extension
$E\colon\B^{\infty}(\ov\rr_+^n)\to\B^{\infty}(\rr^n)$ such
that
  $E\colon \B^{k}(\ov \rr_+^n)\to \B^{k}(\rr^n)$ have bounds depending only on $k$.

\setcounter{thm}{0}\setcounter{equation}{0}
\section{$J$-holomorphic curves and derivatives on curves
}\label{sec3}
In this section, we   first explain how we arrive at the condition  $\|J_i-J_{st}\|<2$
   in  \rt{combvn}. 
Our second result is about $J$-holomorphic curves with parameter. The result is essentially in work of
Nijenhuis-Woolf~\ci{NWsith}. 
 See also Ivashkovich-Rosay~\ci{IRzefo} for another regularity proof and   jets  of $J$-holomorphic
curves. The proof below relies only on some basic facts about the Cauchy-Green operator
and the inversion of $\I+TA\ov{\pd_z}$ discussed in section~\ref{sec2}.
Finally, we will   express partial derivatives through a family of derivatives on  curves.

\smallskip

  Our results are local. Throughout the paper, a real hypersurface $M$ will be
   a relatively open subset  of the boundary of a domain in $\cc^n$, or   a  closed
subset
without boundary in the domain.

Let $\Om$ be a domain in $\rtn$ and $M$ be a (relatively  open) subset of  $\pd\Om$.
  Let $k\geq1$. We say that $X_1,\ldots, X_n$ define an  {\it almost complex structure} $J$ on $\Om$ (resp.~$\Om\cup M$)
of class $\cL C^{k+\all}$,
if $X_j$'s and their conjugates are pointwise $\cc$-linearly independent on $\Om$ (resp.~$\Om\cup M$) and
$X_j$ are of class
$\cL C^{k+\all}$ on $\Om$ (resp.~$\Om\cup M$).
Note that $J_p$ is defined to be the linear map on $T_p\Om$ (resp.~$T_p(\Om\cup M)$)   such that $v+\sqrt{-1}J_pv$ is
in the linear span of   $X_1(p),\ldots, X_n(p)$. The {\it operator norm} of a linear map $A$ from $T_p(\Om)$ into itself is
defined as $\max\{\|Av\|\colon\|v\|=1\}$ with $\|\cdot\|$ being
the euclidean norm on $T_p\Om\equiv\rtn$.
We say that a diffeomorphism $\varphi$ transforms $X_1,\ldots, X_n$
 into $\tilde X_1,\ldots, \tilde X_n$, if $d\varphi( X_j)$ are locally
 in the span of $\tilde X_1,\ldots, \tilde X_n$. 

A linear complex structure $J$ on $\cc^n$ is given by
$$
X_j=\sum_{1\leq k\leq n}(b_{jk}\pd_{\ov z_k}+a_{jk}\pd_{z_k}),\quad 1\leq j\leq n,
$$
where constant matrices $A=(a_{jk})$ and $B=(b_{jk})$ satisfy
$$
 \begin{vmatrix}B & A\\
  \ov A& \ov B\end{vmatrix}
  \neq0.
$$
Denote by $A^t$ the transpose matrix of $A$.
The map     $z=  \ov B^tw+A^t\ov w$ transforms
 $\pd_{\ov w_1}, \ldots, \pd_{\ov w_n}$ into $X_1,\ldots, X_n$, and hence $J_{st}$  into $J$  given by
\eq{abj}
J=
(K^t)^{-1}J_{st}K^t,\quad J_{st}=\begin{pmatrix}0  & \I \\
 -\I &0 \end{pmatrix},\quad
 K=\begin{pmatrix}\RE\, (B+A )  & \IM\,(A+B) \\
 \IM\, (A-B) & \RE\, (B-A ) \end{pmatrix}.
\eeq
 Thus under a
local change of coordinates by shrinking $\Om$ or $\Om\cup M$, an
almost complex structure $J$ is  given by \eq{xjaj} X_j=\pd_{\ov
z_j}+\sum_{1\leq k\leq n} a_{jk}(z)\pd_{z_k}, \quad j=1,\ldots, n
\end{equation}
with operator norm $\|(a_{jk})(z)\|<1$ on $\Om$ (resp. $\Om\cup M$).

\le{jjv}
Let $J^1, J^2$ be two linear complex structures  on $\rr^{2n}$.
 Let $M$ be a  hyperplane in $\rtn$.
There exists  $v\in T_0 M$ such that $J^1v, J^2v$ are in the same connected
component of  $T_0\rtn\setminus T_0M$, provided $T_0M\cap J^1T_0M\neq T_0M\cap J^2T_0M$, or
the  operator norm $\|J^2-J^1\|<2$.
\ele
\begin{proof} To simplify notations, all tangent vectors or spaces are at the origin.
  Let $T(M,J^i)=TM\cap J^iTM$.  
Let $\om_1,\om_2$ be two connected components of $T\rtn\setminus TM$.
Note that  $J^i$ sends one of two
connected components of $TM\setminus T(M, J^i)$
into $\om_1$ and the other     into $\om_2$.
Thus the assertion is trivial, if $T(M, J^1)\neq T(M, J^2)$.  Assume  that they are identical.

 By choosing an orthonormal basis for $T\rtn$,  we  may assume that
  $T(M, J^1)$ is given by $x_n=y_n=0$. Since $M$ contains $x_n=y_n=0$,  then $M$
  is defined by $y_n'=ax_n+by_n=0$ with $a^2+b^2=1$. By a change of orthonormal
  coordinates,
   $M$, $T(M,J^1)$ are defined by $y_n=0$ and  $x_n=y_n=0$ respectively. Write
$$
J^i\binom{\pd_{x'}}{\pd_{y'}}=A_i\binom{\pd_{x'}}{\pd_{y'}}, \quad J^i\binom{\pd_{x_n}}{\pd_{y_n}}
=C_i\binom{\pd_{x'}}{\pd_{y'}}+D_i\binom{\pd_{x_n}}{\pd_{y_n}}.
$$
Here $A_i, C_i,D_i$ are matrices. In particular, $D_i^2=-\I$. We want to show that the coefficients
of $\pd_{y_n}$ in $J^i\pd_{x_n}$ have the same sign. Otherwise,
we can write
$$
D_1=\begin{pmatrix}a_1  & b_1 \\
 -\f{1+a_1^2}{b_1} &-a_1 \end{pmatrix},\quad D_2=\begin{pmatrix}a_2  & -b_2 \\
 \f{1+a_2^2}{b_2} &-a_2 \end{pmatrix}, \quad b_1>0,\quad  b_2>0.
$$
We have $\|D_2-D_1\|<2$. Thus, $b_1+b_2<2$ and $b_1^{-1}+b_2^{-1}<2$,   a contradiction.
 \end{proof}

\begin{exmp}\label{ex+} \rl{jjv} and \rt{combvn} fail easily  for the triplet $\{J_{st}, -J_{st}, \{y_1=0\}\}$.
A less simple  example is in higher dimension. Let $  0\leq t\leq\pi$, and let $J_t$ be defined
  by
$$
X_1^t=(\cos t\, \pd_{x_1}+\sin t\, \pd_{x_2})+i\pd_{y_1},\quad
X_2^t=(-\sin t\, \pd_{x_1}+\cot t\, \pd_{x_2})+i\pd_{y_2}.
$$
 \rl{jjv} and \rt{combvn} fail   for   $\{J_0, J_\pi, \{y_2=0\}\}$ with $\|J_0-J_\pi\|=2$.
Under   new orthonormal coordinates $w_1=(x_2+iy_1)/{\sqrt 2}$, $w_2=(-x_1+iy_2)/{\sqrt 2}$, $J_t$ is given by
\gan
(1+\sin t)\pd_{\ov w_1}-\cos t\pd_{\ov w_2}- (1-\sin t)  \pd_{w_1}-\cos t\,\pd_{w_2},\\
\cos t\,\pd_{\ov w_1}+(1+\sin t)\pd_{\ov w_2}+\cos t\, \pd_{w_1}-(1-\sin t)\pd_{w_2}.
\end{gather*}
The above can be put into   \re{xjaj} with
 $$
 (a_{jk}^t)=
\begin{pmatrix}0  & -\f{\cos t}{1+\sin t}\\ \f{\cos t}{1+\sin t} &0\end{pmatrix}.$$
 Note that $ \|(a_{jk}^t)\|\leq1$. However, we do not know if
 \rt{combvn}  holds for two structures  of the form \re{xjaj} with $\|(a_{jk})\|<1$.
\end{exmp}

  Let $J$ be an almost complex
structure defined by vector fields $X_1,\ldots, X_n$ of class $\cL C^k$
on $\Om$ with $k\geq1$.  A $\mathcal C^1$ map $u\colon\ov\D^+\to\Om$
is called an {\it approximate} $J$-holomorphic curve attached to the curve $u(x,0)$, if
$$
du(\pdoz)=D(z)\cdot X(u(z))+F(z)\cdot \ov{X(u(z))}, \quad |F(z)|=o(|\IM z|^{k-1}).
$$
If $F=0$ and $\ov\D^+$ is replaced by $\D$, $u$ is  called   $J$-holomorphic. Note
that if $f$ is a function on $\Om$ the above equation implies that
$$
\pdoz (f(u(z))=D(z)\cdot( Xf)(u(z))+F(z)\cdot (\ov{X}f)(u(z)).
$$
When $X, J$ are defined by \re{xjaj} and $u$ is  $J$-holomorphic,  the identity becomes
$$
\pdoz (f(u(z))=( Xf)(u(z))\cdot\ov{\pd_zu}.
$$

Next two results deal with the existence of  the two types of curves.
 \le{flat} Let $m\geq0$ be an integer or $m=\infty$. Let $0\leq\all<1$.
Let $J$ be an almost complex structure defined by vector fields
$$X_j=\sum_{1\leq k\leq n} b_{jk}\pd_{\ov z_k}+\sum_{1\leq k\leq n} a_{jk}\pd_{z_k}, \quad j=1,\ldots, n.$$
Assume that $A=(a_{jk}), B=(b_{jk})$ are of class $\B^{m+\all}(\Om)$.
 Assume that   $u_0\colon(-1,1)\times P\to K$ is  of class
 $\hB^{l+1+\all,j}((-1,1),P)$, and  $K$ is a compact subset  of $\Om$.
 Let  $l\geq 0, j\geq0$,
  $j+l\leq m$, and $0<r<1$.   There
 exists a 
 map  $u\colon[-r,r]\times[-\del,\del]\times P\to\Om$
of class $\hB^{l+1+\all,j}([-r,r]\times[-\del,\del],P)$ satisfying the following.
 \bppp
 \item  $u(x,0,t)=u_0(x,t)$ and
\ga
\label{duce}
du(\pdoz)=D(z,t)\cdot X(u(z,t))+F(z,t)\cdot\ov{X(u(z,t))},\\
\label{duce+}
|F(z,t)|=o(|y|^{l}), \  \alpha=0; \quad |F(z,t)|=O(|y|^{l+\all}),\  0<\alpha<1.
\end{gather}
 \item Let $e_l>0$.
On  $[-r,r]\times[-\del,\del]\times P$ the norms of $u, D,  F$  satisfy     $$\|u\|_{l+1+\all,j}
+\|(D,F)\|_{l+\all,j}\leq C_l^*\max(\|u_0\|_{l+1+\all,j},\|u_0\|_{l+1+\all,j}^{l+j+2}) +e_l.$$
 \eppp Moreover, $ C_0^*(1+\|u_0\|_{1,0})\del>1$ and $C_l^*$ depends on $K,
 \Om$ and
$$
  |A|_{l+\all},\quad |B|_{l+\all},\quad \inf_{\Om}  {\small \begin{vmatrix}B \ \ A\\
  \ov A\ \ \ov B\end{vmatrix}}.
  $$
\ele
\begin{proof} We suppress the parameter $t$   in all expressions.
We first determine a unique set of   coefficients $a_1(x),\ldots,
a_{l+1}(x)$  such that as a power series  in $y$,
$u(x,y)=u_0(x)+\sum_{i=1}^{l+1}a_i(x)y^i$ satisfies
\re{duce}-\re{duce+}. It is convenient to regard $u$
as the real map $(x,y)\to (\RE u,\IM u)$, which is still denoted by $u$,
 and  rewrite the equations as
\eq{dujud}
du(\pd_y)=J(u)(du(\pd_x))+   F_1(x,y)\cdot \pd_*,\quad F_1(x,y)=o(|y|^l).
\eeq
Here $\pd_*=(\pd_{u_1},\ldots,\pd_{u_{2n}})$ is evaluated at $u(z,t)$.
In the matrix form, let $J$ be the matrix defined by \re{abj}.
Then we need to solve
$$
\pd_yu=\pd_xu J(u)+F_2(x,y), \quad |F_2(x,y)|=o(|y|^l).
$$
We  solve the equation formally, which determines
$a_1(x),\ldots, a_{l+1}(x)$ uniquely, and then apply the Whitney extension (\rl{whit}).
 This gives us a map $u$ from $([-r,r]\times[-1,1])\times P$ into $\rr^n$ of class
  $\cL C^{l+1+\all,j}$  satisfying  the stated
  norm estimate.   By $|u(x,y)-u(x,0)|\leq C_0^*
  (\|u_0\|_{1,0}+e_1)|y|$ and the compactness of $K$, we find
  $\del>0$ such that    $u$ maps $[-r,r]\times [-\del,\del]\times P$ into $\Om$.
We have obtained \re{dujud}. Thus
$$
2 du(\pdoz)=du(\pd_x)+idu(\pd_y)=du(\pd_x)+iJ(u)(du(\pd_x))+iF_1(x,y)\pd_*.
$$
Note that  $du(\pd_x)+iJ(u)(du(\pd_x))=  D_1(z)\cdot X(u(z))$.
Write $\pd_*$ in terms  of
 $X_i,\ov X_j$  by using the inverse of
$ {\tiny\begin{pmatrix}B \  A\\
  \ov A\ \ov B\end{pmatrix}}$.
We get \re{duce}-\re{duce+}. We can estimate  the norms of $D,F$ via
$  D_1, F_i$ and the inverse matrix.
\end{proof}

We have defined $\hB^{k+\all,j}$, $\|\cdot\|_{k+\all,j}$ in section~\ref{sec2}.
Following~\ci{NWsith}, we define for $j\leq k$
\gan
\B^{k+\all, j}(\ov\Om, P)=\bigcap_{0\leq l\leq j}\hB^{k-l+\all,l}(\ov\Om, P),\quad
{\mathbf|}  u|_{k+\all,j }=\max_{0\leq l\leq j} \|u  \|_{k-l+\all,l}.
\end{gather*}
 One can see that $ \B^{k+\all, j}(\ov\Om, P)$ is complete.
By   assumptions on $\Om,P$, we see  that
  $$\B^{k+\all,k}(\ov\Om,P)\supset\cL C^{k+\all}(\ov\Om\times   P).
$$
In particular, if $f\in\cL C^{k+\all,j}\cap\cL C^1$ and $u\in\B^{k+\all,j}(\ov\Om,  P)$,
then $f\circ u\in\B^{k+\all,j}(\ov\Om,  P)$ whenever  the composition is well-defined.
In general, let
 $\varphi(x,t)=(\tilde\varphi(x,t),t)$ with $\tilde\varphi$ being a map
 from $\Om\times P$ into $\Om'$  of class $\B^{k+\all,j}\cap\mathcal C^1$.  Then
$$
|v\circ\varphi|_{k+\all,j} \leq C(1+|\tilde\var|_{1,0}+|\tilde\var|_{k+\all,j} )^{1+k+j}|v|_{k+\all,j}.
$$

Let $\D$ be the unit disc in $\cc$, $\D_r$ the disc of radius $r$, and
$\D_r^+=\D_r\cap\{\IM z>0\}$.
The following result gives coordinate  maps  in $J$-holomorphic curves.
\pr{corZ} Let $0<\all<1$ and let $j\geq1$ be an integer. Let $J$ be an almost
complex structure defined by vector fields
$$X_i=\sum_{1\leq k\leq n} b_{ik}\pd_{\ov z_k}+\sum_{1\leq k\leq n} a_{ik}\pd_{z_k},
 \quad j=1,\ldots, n.$$
Assume that $A=(a_{ik}), B=(b_{ik})$ are of class $\B^{j+1+\all}(\Om)$.   Let
$M\subset\Om$ be  a real hypersurface   of class $
\cL C^{j+2+\all}$.  Let $e\colon M\to\cc^n$ be a $\cL C^{j}$ map such that
 $e\cdot X=e_1X_1+\cdots+e_nX_n$ is
   not tangent to $M$ at each point of $M$. Let $0\in M$.
  There exist
  two $\cL C^j$ diffeomorphisms   $u, R$  from $\D_r^n$ into $\Om$
    satisfying
the following.
\bppp
\item For each $t\in\D_r^{n-1}$,
  $u(\cdot,t)$
is $J$-holomorphic and embeds  $ \D_r$ onto  $D (t)$.
\item  $u(0,t)$ is in $M$, and  $D(t)$ intersects
$ M$ transversally along a curve $\gamma(t)$.  Also, $u(0)=0$ and
 $du(0,t)(\pd_{\ov\zeta})=(e\cdot X)(u(0,t))$.
\item
    $R(\cdot,t)$ sends $\D_r^+, (-r,r),
 \D_r$ into $\Om^+\cap D(t), M\cap D(t), D(t)$, respectively.
 And $R(0)=0$.
 \item  If $A\in\cL C^{k+\all}$,
 $M\in\cL C^{k+1+\all}$ and $k> j$,   the 
 $u$,  $ R$ are in $\B^{k+1+\all,j}(\D_r,\D_r^{n-1})$.
\eppp
Here $r$  depends only on $\inf_{\Om} {\tiny\begin{vmatrix}B \ A\\
  \ov A\ \ov B\end{vmatrix}}$,  $M$, $e$,
$j$, $\all$, $|(A,B)|_{j+\all}$,
 and the diameter of $\Om$.
\epr
 \begin{proof} 
Introducing the new coordinates $w$ by
$z=  \ov B^t(0)w+A^t(0)\ov w$, we may assume that  $A(0)=0$
  and $B=\I$. Thus we obtain
$$ 
   X_j=\pd_{\ov
z_j}+\sum_{1\leq k\leq n} a_{jk}(z)\pd_{z_k}, \quad j=1,\ldots, n.
$$ 
Applying a unitary change of coordinates, we may assume that
$T_0M$ is given by
 $y_n=0$.
 By a
 change of coordinates which is tangent to the identity and
 of class $\cL C^{k+1+\all}$, we may assume
 that $M$ is in $y_n=0$.  By dilation,  $\Om=\D_2^n$ and on it we have $\|A(z)\|<1/4$ and
 $ |A |_{j+1+\all}<1/{C_{*}}$. Here $C_{*}$ will be determined. Finally, by a dilation in
 $\D$, we achieve $\|e(x)\|<1/{4}$ on $M$.

\smallskip

\noindent {\it Existence in $\B^{j+1+\all,j}$  class.}
We first find   $u$.
Recall that $u\colon \D\to\Om$ is $J$-holomorphic, if $du(\pd_{\ov\zeta})$ is
in the span of $X_i's$. Then the equations  are
$$
\pd_{\ov\zeta}u_i=\sum_{1\leq l\leq n}a_{li}(u)\ov{\pd_{\zeta} u_l}, \quad i=1,\ldots, n.
$$
In column vectors, they become
\eq{jeqa}
\pd_{\ov \zeta}u=A^t(u)\ov{\pd_{\zeta} u}.
\eeq
 At the origin,
 $X_j(0)=\pd_{\ov z^j}$.
  Let $\tilde e_1,\ldots,\tilde e_{n-1}$ be  the standard base of $\cc^{n-1}\times0$.
Then $\tilde e_1,\ldots, \tilde e_{n-1}$,
 $e(0,t)$
 are  $\cc$-linearly independent.
For  $t\in P=\D^{n-1}$, we look for a $J$-holomorphic curve $u$ satisfying
$
u=\Psi(u)$ with
\eq{Psiu}  \Psi(u)(\zeta)=t\cdot(\tilde e_1,\ldots,\tilde e_{n-1})/n+\zeta \ov {e(0,t)}
+\Phi(u)-P_1\Phi(u).
\eeq
Here $\Phi(u)=
T_{\D}(A^t(u )\ov{\pd_\zeta u})$ and
 $P_1\Phi(u)(\zeta,t)=\Phi(u)(0,t)+\pd_{\zeta}\Phi(u)(0,t)\zeta$.
Let
$\mathcal B_1$ be the  closed  unit ball
in
$\mathcal B=[\B ^{j+1+\all,j}
  (\ov\D,  P)]^n
$
equipped with norm
$
|u|_{j+1+\all,j}.
$ When $u\in\mathcal B_1$, $\Phi(u)$ is in $\mathcal B$ (see \rl{tspa}).
   Then $P_1\Phi(u)$ is continuous and
 of class $\cL C^{j}$ in $t$. It is also a polynomial
 in $\zeta$. In particular, $P_1\Phi(u)$ and  $\Psi(u)$ are in $\mathcal B$.
One can verify that if
$|A|_{ j+1+ \all}<1/C_{*}$ on $\D_2\times P$, then $u\to \Psi(u)$ is a contraction map from
$\mathcal B_1$ into itself (here we need $A$ to be in $\cL C^{j+1+\all}$ instead of $\cL C^{j+\all}$). We take $u\in\mathcal B_1$ to be its fixed point.

Recall that after dilation,  $|A|_{j+1+\all}<1/{C_{*}}$ and $j\geq1$. Then
$\Phi(u)$ and $P_1(\Phi(u))$ have
small $\cL C^1$ norms on $\cL B_1$ in $\zeta,t$ such that $u$   is a $\cL C^1$
diffeomorphism in $\zeta,t$.

\smallskip

\noindent
{\it Higher order derivatives.}
We have obtained a solution  $u\in\B^{j+1+\all,j}(\ov\D,P)$.
Assume now that  $A\in\B^{k+\all}(\D_2^n)$ with $k>j$.
We want to show   a stronger result:
Assume that for all $l\leq j$,
$\pd_t^lu(\zeta,t)$ are continuous on $\D\times   P$
and distributional derivatives
$\pd_\zeta\pd_t^lu(\cdot,t)$ have bounded   $L^p(\D)$ norms
on $ P$ with    $p>2$.
Assume that $u(\cdot,t)$ is $J$-holomorphic on $\D$. Then $u\in\B^{k+1+\beta,j}(\D_r,
P)$
for $r<1$ and  $ \beta=\min(\all, 1-{2}/{p})$.

Indeed, \re{jeqa} implies that the  first-order derivatives of $\pd_t^lu(\cdot,t)$
    have bounded $L^p(\D)$ norms on $P$. By Morrey's inequalities,
$u\in\hB^{\beta,j}(\D_r,P)$ for any $r<1$. (See Lemma 7.16 and Theorem 7.17 in~\ci{GTzeon}, pp.~162-163.)

Fix $\zeta_0\in\D$. Let
$u=\tilde u+A^t(u(\zeta_0,t))\ov{\tilde u}$ and $\tilde u(\zeta_*,t)=u_*(\zeta,t)$ with
$\zeta_*=\zeta_0+\mu\zeta$. Here   $0<\mu<\yt(1-|\zeta_0|)$ will be determined. We get
on $\D$
\ga\label{jeqa+}
\pd_{\ov \zeta}u_*=A_*^t(\zeta,t)\ov{\pd_{\zeta} u_*},\quad A_*(0,t)=0,\\
A_*(\zeta,t)=[A(u(\zeta_*,t))-A(u(\zeta_0,t))][\I-\ov{A(u(\zeta_0,t))}
A(u(\zeta_*,t))]^{-1}.\label{jeqa++}
\end{gather}
Let $\chi$ be a smooth function with  support in $\D_{1/4}$.
 Let $v=\chi u_*$.
Multiply \re{jeqa+} by $\chi$ and rewrite it as
\eq{deqv}
\pd_{\ov\zeta} v-A_*^t(\zeta,t)\ov{\pd_{\zeta}v}=u_*\pd_{\ov\zeta}\chi
-A_*^t(\zeta,t)\ov{ u_*\pd_{\zeta}\chi}.
\eeq
Let $\tilde\chi$ be a smooth function  with compact support in $\D$. We  also assume that
 $\tilde\chi=1$ on $\D_{1/4}$  and $|\tilde\chi|_1<5$. Replacing
$A_*$ by $\tilde\chi A_*$, we may assume that
$A_*(\cdot,t)$ has compact support in $\D$.
Using \re{jeqa++}, we get  for $\zeta,\zeta'\in\D$,
\gan
|A_*(\zeta,t)|\leq C
|A(u(\cdot,t))|_{\beta}\mu^{\beta},\\
|A_*(\zeta',t)-A_*(\zeta,t)|\leq C
|A(u(\cdot,t))|_{\beta}\mu^{\beta}|\zeta'-\zeta|^{\beta}.
\end{gather*}
Therefore, $\|A_*\|_{\beta,0}\leq
C
|A\circ u|_{\beta,0}\mu^{\beta}<\e_{\beta}$. Here
$\e_{\beta}$ is the constant
in \rl{isa} and $\mu$ is sufficiently small. Apply $T=T_\D$ to \re{deqv}. Since $v$ has compact support, then
$$
v-T(A_*^t\ov{\pd_\zeta v})= T(u_*\pd_{\ov\zeta}\chi -A_*^t\ov{ u_*\pd_{\zeta}\chi}).
$$
Write the right-hand side as $w$ and solve for $v=(\I-TA_*^t\ov{\pd_{ \zeta}})^{-1}w$. Since $u_*$  is in
  $\hB^{\beta,j}(\ov\D,P)$, then $A_*, w$
 are in $ \hB^{\beta,j}(\ov\D,P)$. By \rl{isa}, $v$ and hence $u$ are in $\hB^{1+\beta,j}$.
Repeating the procedure, we get $u\in\hB^{k+1-j+\beta,j}$. Also $u\in\hB^{k+1-l+\beta,l}$
for all $l\leq j$.
This shows that $u\in\B^{k+1+\beta,j}$.

\smallskip

\noindent {\em End of the proof}.  We assume that $\Om=\D_2^n$ and that  $\Om^+, M$
are subsets defined by $y_n>0$ and $y_n=0$, respectively. Let $\ov {e(0,t)}=(a, b'+ib'')$.
Since $ e(0,t)\cdot\pd_{\ov\zeta}$ is not tangent to $M$, then $b'+ib''\neq0$. Without loss of
generality,
we may assume that $b'\geq|b''|$.
We have $u=\Psi(u)$. By \re{Psiu}, $D(t)\cap M$ is defined by
\eq{bxbe}
b''\xi+b'\eta=F(\xi,\eta,t), \quad F(\xi,\eta,t)=
\IM\{P_1\Phi_n(u(\xi+i\eta,t))-\Phi_n(u(\xi+i\eta,t))\}.
\eeq
We already know that $F\in\B^{k+1+\all,j}(\D_r,P)$. We  may also achieve $|\pd_\eta F|<b'/2$, by
assuming $|A|_{j+\all}<1/C_*
$.
By the implicit function theorem, \re{bxbe} has a   solution
 $
 \eta=h(\xi,t)$ for $  |\xi|<r/c,  t\in P$.
 Now $$(\pd_\xi h,\pd_th)=(b'-\pd_\eta F(\xi,\eta,t))^{-1}\bigl(\pd_\xi F-b'',\pd_t F\bigr)$$
 implies that $\pd_t^l
 h\in\cL C^{k+1+\all-l}$ for all $l\leq j.$
  On $\D_{r/c}\times P$, define
 $$
 R(\zeta,t)=u(\xi+i(\eta+h(\xi,t)), t).
 $$
Then $R(\cdot,t)$ sends $\D_{r/c}^+$ into $D^+(t)$. Replace $R(\zeta,t)$ by $R(\zeta/c,t)$.
The remaining assertions can be verified
easily.
\end{proof}

We remark that the above $R$ is
not  $J$-holomorphic.

\medskip

It is well-known that via the Fourier transform,   the boundedness of derivatives of
 a function on all lines parallel to coordinates axes
yields some smoothness of the function  in all variables  (see Rudin~\ci{Runion}, p.~203).
To limit the loss of derivatives, we will use
the  Fourier transform only on curves. This requires us to bound  derivatives
of a function on a larger family of
curves.

Let  $\gaa$ be a $\cL C^k$ curve in $\rr^n$, and let $f$ be a function
 of class $\cL C^k$ on $\rr^n$. We have
\eq{dgk}
\pd_{t}^kf(\gaa(t))= (\gaa'(t)\cdot\pd)^k f) (\gaa(t))
+\sum_{1\leq |\beta|<k}\!
Q_{k,\beta}(  \pd_t^{(k+1-|\beta|)}\gaa)(\pd^\beta f)(\gaa(t)).
\eeq
Here   $Q_{k,\beta}$ are polynomials, $\pd^{(k)}$ denotes
derivatives of order $\leq k$, and
 $$
 v\cdot \pd=v_1\pd_{x_1}+\cdots+v_n\pd_{x_n}.
 $$
 \le{det} Let $k$ be a positive
integer and let $\e>0$. \bppp
\item
There exist $N$ vectors $v_j=(1,v_j')\in\rr^n$
such that $|v_j'|<\e$ and
\eq{xiac}
\pd^\all=c_{\all,1}(v_1\cdot\pd)^k +\cdots+ c_{\all,N}(v_N\cdot\pd)^k, \quad  |\all|=k.
\eeq
\item If $v_1,\ldots, v_N$ satisfy \rea{xiac},
there exists $\del>0$ such that if $|u-v|<\del$, then
\eq{xiac+}
\pd^\all=Q_{\all,1}(u)(u_1\cdot\pd)^k +\cdots+ Q_{\all,N}(u)(u_N\cdot\pd)^k, \quad  |\all|=k.
\eeq
\eppp
Here $Q_{\all,j}$ are rational functions
with $Q_{\all,j}(v)=c_{\all,j}$. And $N$ depends only on $k,n$.
\ele
\begin{proof} (i). Equivalently, we
 need to verify \re{xiac}-\re{xiac+} when $\pd$ is replaced by $\xi\in\rr^n$.
It holds for $n=1$. Assume that it holds when $n$ is replaced by $n-1$.
For $\xi_n^k$, we take   distinct non-zero constants $\la_1,\ldots, \la_{k}$.
Then $\xi_n^k$ is in the linear span of
$\xi_1^k$, $(\xi_1+\la_1\xi_n)^k,\ldots, (\xi_1+\la_k\xi_n)^k$.
 Let $\xi_n^jP(\xi_1,\ldots, \xi_{n-1})$ be a
 monomial of degree $k>j$. Then by the induction assumption
$$
\xi_n^jP(\xi_1,\ldots, \xi_{n-1})= \xi_n^j [c_1 (v_1\cdot \xi)^{k-j}+\cdots
+c_l (v_l\cdot\xi)^{k-j}].
$$
Here $v_j=(1,v_j'',0)$ with $|v_j''|<\e/2$.
Then  $\xi_n^i(v_l\cdot \xi)^{k-i}$  are in the linear span of
$(v_l\cdot\xi)^k$,  $  (v_l \cdot \xi+\la_j \xi_n)^k$  with $ j=1,\ldots, k.
$
Note that $\la_j$ can be arbitrarily small. Thus, (i) is verified.

(ii). For $|\all|=k$ we have expansions
\aln
\xi^\all&=\sum_{1\leq j\leq N}c_{\all,j}(v_j\cdot\xi)^k,\quad \xi^\all
=\sum_{1\leq j\leq N}c_{\all,j}(u_j\cdot\xi)^k+\sum_{|\beta|=k}\widetilde Q_{\all \beta}(v-u)\xi^\beta.
\end{align*}
Clearly, $\widetilde Q_{\all \beta}(0)=0$. Moving the last sum to the left-hand side
and inverting
$\I-(\widetilde Q_{\all \beta})$ yields \re{xiac+}.
\end{proof}

We now use \re{dgk} to   estimate partial derivatives via   derivatives on curves.  Set $t'=(t_2,\ldots, t_n)$
and $t=(t_1,t')$.
\pr{ftva} Let $k, N$ be   positive integers.
For $1\leq j\leq N$, let $R_j$ be  $\cL C^{1}$ diffeomorphisms from $\Om_j\subset\rr^n$
 onto an open subset $\Om$ of $\rr^n$. Assume that $
 R_j(\cdot,t')\in\cL C^{k}$ and $R_j(0)=0$. Suppose that  at $0\in\Om$
\eq{pda0}
\pd^\all=\sum_{1\leq j\leq N}c_{\all,j}(\pd_{t_1}R_j(0)\cdot\pd)^{|\all|}, \quad 1\leq |\all|\leq k.
\eeq
Let $f\in \cL C^0(\Om)$. Then the following hold.
\bppp
\item Let $f$ be of class $\cL C^k$ near $0\in\Om$.
There exist rational functions $Q_{\all,i,j}$ such that for
$x=R_j(t^j)$ near $0$ and   $|\all|=m\leq k$ with $m\geq1$, \eq{pda}
\pd^\all f(x) =\sum_{i=1}^{m}\sum_{j=1}^NQ_{\all,i,j}
\Bigl(\pd_{t_1^1}^{(m-i+1 )} R(t^1),\ldots,\pd_{t_1^N}^{(m-i+1
)}R(t^N)\Bigr)\pd_{ t^j_1}^if(R_j(t^j)). \eeq
\item
Suppose that $R_j$ are affine, i.e.  $R_j(t)-R_j(y)=R_j(t-y)$
wherever they are defined.
 Suppose that   $L_{t_1}^\infty$ norms of one-dimensional distributions
 $\pd_{t_1}^m(f\circ R_j)(\cdot,
t')$  are  bounded in $t'$ for all $m\leq k$.
Then 
near $0$, $\pd^\all f$ are
Lipschitz functions
for all $|\all|<k$.
\item  Let $R_j$ be of class $\cL C^{k+1}$ near $0\in\rr^n$ and
let $n<p<\infty$.
Suppose that
    $L_{t_1}^p$ norms of one-dimensional distributions
    $\pd_{t_1}^m(f\circ R_j)(\cdot,
t')$  are bounded  in $t'$
  for all $m\leq k$.
Then near $0$,
$f$ is of class $ \cL C^{k-\f{n}{p}}$. 
\eppp
\epr
\begin{proof}
(i) follows from \re{dgk} and \re{xiac+}, by hypothesis \re{pda0}.

(ii).
 Applying  dilation and replacing $f$ by $\chi f$, we may assume that
$f$ has compact support in $\Del^n$.
Let $\chi_\e(x)=\e^{-n}\chi(\e^{-1}x)$ for a smooth function $\chi$ with support
in $\Del^n$ and  $\int\chi\, dx=1$. Let $f_\e(x)=\int f(y)\chi_\e(x-y)\, dy$ and
$f_{\e,j}=f_\e\circ R_j$.

Changing variables via $R_j$,
we get
$$
f_{\e,j}(t )=\int f(R_j(t)-R_j(y))\chi_\e(R_j(y))\det R_j'(y)\, dy.
$$
Using $R_j(t)-R_j(y)=R_j(t-y)$,  we get   $|f_{\e,j}(\cdot,
t')|_{k}<C$ for $C$ independent of $\e$ and $t'$.
In \re{pda}, we substitute $f_\e$ for $f$.
Therefore, $\pd^\all f_\e$ are bounded near $0$. We can find a sequence $f_{\e_j}$
such that as $\e_j$ tends to $0$,  $\pd^\all f_{\e_j}$
converges uniformly for $|\all|<k$, and the Lipschitz norms of
$\pd^\all f_{\e_j}$ are bounded by
a constant.   Since $f_\e$ converges to $f$ uniformly as $\e\to0^+$ then
 $\pd^{k-1}f\in Lip_{\, loc}$.

 (iii).  For the $f$, we
define a distribution $T_jf$ by
 $$
 T_jf(\phi)=(-1)^k\int_{\rr^n} f\circ R_j(t)\pd_{t_1}^k(\phi(R_j(t)))\, dt.
 $$
 Here $\phi$ are test functions  supported in $\Del^n_\e$ with $\e$ small.
 It is clear that  defined near $0$, $T_jf$ is a distribution of order $(\leq) k$.
 Integrating in $t_1$-variable first and throwing  the one-dimensional  derivative onto
  $f\circ R_j$ yields
 \aln
 |T_jf(\phi)|&\leq C \int_{\rr^{n-1}}\|\pd_{t_1}^k[f\circ R_j](\cdot, t')\|_{L_{t_1}^p}\|\phi\circ R_j(\cdot,t')
 \|_{L_{t_1}^q}\, dt'\\
&\leq  C_1\int \|\phi\circ R_j(\cdot,t')
 \|_{L_{t_1}^q}\, dt'\leq C_2\|\phi\circ R_j\|_{L^q}\leq  C_3\|\phi\|_{L^q}.
 \end{align*}
Here the second last inequality is obtained from the H\"older inequality
and $\supp\phi\subset \Del^n_\e$.
 Hence near $0$,  $T_jf\in L^{p}$ when $p>1$.
 Next we find  a differential operator $P_{j,k}(\pd)$ of
 order $k$ such that $P_{j,k}(\pd)f=T_jf$.   To find it,
 we use a smooth function $g$ to obtain
   \aln
T_jg(\phi)&=
 \int \pd_{t_1}^k[g\circ R_j(t)]\phi(R_j(t))\, dt
 =\int (\phi\tilde P_{j,k}(\pd )g)\circ R_j(t)\, dt\\
 &=\int [\det((R_j^{-1})')\tilde P_{j,k}(\pd )g]\phi\, dx
 \df  (P_{j,k}(\pd) g)(\phi).
 \end{align*}
 Since $R_j\in \cL C^{k+1}$, it is easy to see that
 \gan
 P_{j,k}(\pd)=\det((R_j^{-1})')\tilde P_{j,k}(\pd )=\sum_{|\all |\leq k}a_{j,k,\all}\pd^\all,
 \quad a_{j,k,\all}\in \cL C^{|\all|}.
 \end{gather*}
 The last assertion implies that $P_{j,k}(\pd)$ has order $k$.  The definition of $T_jf$ and  identity
    $P_{j,k}(\pd) g=T_jg$ implies that as distributions defined  near $0$,
  $P_{j,k}(\pd)f=T_jf$.

   Note that
 \gan
 \sum_{|\all |=k}a_{j,k,\all}(x)\pd_x^\all=\sum_{|\all|=k}C_\all\det((R_j^{-1})')(\pd_{t^j_1}R_j(t^j))^\all \pd_x^\all, \quad C_\all\neq0.
 \end{gather*}
 Here $t^j=R_j^{-1}(x)$.
 Combining with \re{pda0}, we get for $g\in\cL C^k$ and $1\leq|\all|\leq k$,
  $$
\pd^\all g=\sum_{1\leq i\leq m}\sum_{1\leq j\leq N}b_{\all, j,i}P_{j,i}(\pd)g,\quad
b_{\all,j,i}\in\cL C^{i}.
 $$
 The last assertion, combined with $\ord P_{j,i}(\pd)\leq i$, $a_{j,k,\all}\in\cL C^{|\all|}$ and $P_{j,i}(\pd)f\in L^p$,
  implies
   that near $0$, $\pd^\all f$ are in $ L^p$ for
  $1\leq |\all|\leq k$;
  by a Sobolev embedding theorem (\cite{Honize}, p.~123), $f\in\cL C^{k-1+\beta}$
  with $\beta=1-\f{n}{p}$.
\end{proof}

\setcounter{thm}{0}\setcounter{equation}{0}
\section{Cauchy-Green operator on domains with parameter
}\label{sec3+}

The following result is certainly  classical; see~\ci{Vesitw}, section 8.1 (pp.~56-61).
For the convenience of the reader, we present details for  a parameter version.
 Recall that
$P$ is the closure of a  bounded open set in a euclidean space and
      two points $a,b$ in $  P$ can be connected
by a smooth curve in $   P$ of length at most $C|b-a|$.

\le{ctva}
Let $\tau$ be a complex-valued function on $\ov\D^+\times P$   of class $
  \B^{k+1+\all,0}(\ov\D^+,P)$. Suppose that for $z,z'\in\D^+$ and $t\in P$,
 \eq{rz't}
|\tau(z',t)-\tau(z,t)|\geq |z'-z|/C.
 \eeq
\bppp\item Let $f$ be a continuous function on $ [-1,1]\times P$. Let
$$
  C_0f(z,t)=\f{1}{2\pi i}\int_{-1}^1\f{f(s,t)}{\tau (s,t)
 -\tau (z,t)}  \, ds,\quad z\in\D^+.
 $$
Then $|\pd_z^kC_0f(z,t) |\leq C_k|f|_{0}/{|\IM z|^{k+1}}$, where $c_k$
depends only on $\|\ta\|_{k,0}$.
\item If
  $f$ is a function  of class $\B^{k+\all,0}([-1,1],P)$,
then $C_0f$ extends continuously to   $(\D\cup(-1,1))\times P$.
Moreover,   $C_0f\in\B^{k+\all,0}(\ov\D_r^+,P)$  for   $r<1$ with $|C_0f|_{k+\all,0}\leq C|f|_{\all,0}$.
\item
 Let $f$ be a function  of class $\B^{k+\all,0}(\ov\D^+,P)$. For $ z\in\D^+$, define
\gan
  S_0f(z,t)=-\f{1}{\pi}\lim_{\e\to0}\int_{ \{\zeta\in\D^+\colon
  |\ta(\zeta,t)-\ta(z,t)|>\e\}}\f{f(\zeta,t) }{(\tau (\zeta,t)
 -\tau (z,t))^2}  \, d\xi d\eta,\\
  T_0f(z,t)=-\f{1}{\pi}\int_{\D^+}\f{f(\zeta,t)}{\tau (\zeta,t)
 -\tau (z,t)}  \, d\xi d\eta.
 \end{gather*}
Then     $S_0f\in\B^{k+\all,0}(\ov\D_r^+,P)$
and $T_0f\in\B^{k+1+\all,0}(\ov\D_r^+,P)$ for $r<1$ with $|S_0f|_{k+\all,0}+|T_0f|_{k+1+\all,0}\leq C|f|_{k+\all,0}$.
\eppp
\ele
\begin{proof} (i).
Note that \re{rz't} implies that $|\ta(z,t)-\ta(s,t)|\geq \IM z/C$
for $-1\leq s\leq1$ and $z\in\D^+$.
The proof is straightforward by taking derivatives
in $z,\ov z$ directly onto the kernel.

 (ii). Let $z=x+iy$. 
 Let $\chi$ be a smooth function
with compact support in $(-1,1)$.
Replacing $f(x,t)$ with $\chi (x)f(x,t)/{\pd_x\tau (x,t)}$,
  it suffices to get the  norm estimate   on $\D_r\times P$
 for
\al\label{tcfg}
 C_0f(z,t)&=\f{1}{2\pi i}\int_{\pd\D^+}
 \f{f(\zeta,t)}{\tau (\zeta,t)
 -\tau (z,t)}  \, d\tau (\zeta,t)\\
\nonumber &=\f{1}{2\pi i}\int_{\pd\D^+}
 \f{f(\zeta,t)-f(x,t)}{\tau (\zeta,t)
 -\tau (z,t)}  \, d\tau (\zeta,t)+\e   f(x,t).
\end{align}
Here the differentiation and integration are in $\zeta$.
And $\e=1$, if  $\tau(\cdot,t)$ preserves the orientation
 of $ \D^+$; otherwise $\e=-1$. From \re{rz't} and $\ta\in\B^{1,0}(\D^+,P)$,
 we know that $\e$ is independent of $t$.
  Let $C_1f$ denote the second integral in \re{tcfg}.
We   denote  a $j$-th derivative  in $x,y$ by $\pd^j$.
 In what follows,
the norms   $|\cdot|_{j+\all,0}$ for $f,\tau $   are on $\D^+$, and
norms   $|\cdot|_{j+\all,0}$ for $ C_0f$ are on $\D_r^+$ with $r<1$. These norms  will be
 denoted
by the same notation $|\cdot|_{j+\all}$.
Since   $t$ is fixed, we suppress it in all expressions.
All constants are independent of   $t$.

 That $C_1f$ extends continuously to $\ov\D^+\times P$
follows   from the continuity of $f$ and
$$
\left|\f{f(s)-f(x)}{\tau (s)
 -\tau (z)}\right|\leq C|f|_{\all}|x-s|^{\all-1}.
$$
Take \re{tcfg} as the definition of $C_0$. Differentiating  it   gives
\eq{pdzg}
\pd C_0f(z )=\f{\pd \tau (z )}{2\pi i}\int_{\pd\D^+}\f{f(\zeta )-f(x )}{(\tau (\zeta )
 -\tau (z ))^2}  \, d\tau (\zeta ).
\end{equation}
Using $|\tau (s )-\tau (z )|\geq (|s-x|+|y|)/C$, we get
$$
|\pd C_0f(z )|\leq\|\tau \|_{1}
\int_{\pd\D^+}\f{C|f|_\all|s-x|^\all}{|y|^2+|s-x|^2}\, ds
\leq C_\all'|f|_\all|y|^{\all-1}.
 $$
By a type
of Hardy-Littlewood lemma, we obtain
$|C_0f|_{\all,0}\leq C|f|_{\all,0}$.
For higher derivatives of $C_0f$, we differentiate  \re{tcfg}
in $z$ variable and transport derivatives to $f$ via  integration by parts.
We get  for $|I|=k$
\eq{pdzg1}
\pd^I C_0f(z )=\sum_{1\leq |J|\leq |I|}\f{\pd^J\tau (z )}{2\pi i} \int_{\pd\D^+}
 \f{f_{IJ}(\zeta )}{\tau (\zeta )
 -\tau (z )}  \, d \tau (\zeta ).
\eeq
Here $f_{IJ}(s )$  are polynomials in $(\pd_s\tau (s ))^{-1},
\pd_s^lf(s ), \pd_s^{l+1}\tau (s )$ with $l\leq k$. As before, we have the continuity of
$$\f{1}{2\pi i}\int_{\pd\D^+}
 \f{f_{IJ}(\zeta )}{\tau (\zeta )
 -\tau (z )}  \, d \tau (\zeta )=\f{1}{2\pi i}\int_{\pd\D^+}
 \f{f_{IJ}(\zeta )-f_{IJ}(x )}{\tau (\zeta )
 -\tau (z )}  \, d\tau (\zeta )+\e   f_{IJ}(x ).
 $$
 Differentiating the integral in \re{pdzg1} one more time
we get a formula analogous
to \re{pdzg}.
   As in   case $k=0$, we can verify that the $\cL C^\all$
norms of  $\pd^kC_0f(\cdot, t)$ on $\ov\D_r$
are bounded. 

(iii).
We first show that $S_0f\in\B^{\all,0}(\ov\D^+,P)$.

 Let $ -2i d\xi\wedge d\eta=A(\zeta,t)\, d\tau(\zeta,t)\wedge d\ov{\tau(\zeta,t)}$. Let $\chi$ be
a smooth function with compact support in $\D^+_{r'}\cup(-r',r')$, $0<r<r'<1$.
 Replace  $f(\zeta,t)$ by
 $\chi(\zeta) f(\zeta,t)A(\zeta,t)$.
 We may
  reduce to the case that $f(\cdot,t)$ is supported in $\ov\D^+_{r'}$
  with $r<r'<1$. We may also replace the domain of integration by a smooth domain $D$ with
  $\D_r^+\subset \ov D\subset\ov\D^+_{r'}$.
 Again, we suppress the parameter $t$ in all expressions and write
\eq{gztf}
S_0f(z )=\f{1}{ 2\pi i } \int_{D}\f{(f(\zeta )-f(z ))\, d\tau(\zeta )
\wedge d\ov{\tau(\zeta ) }}{(\tau (\zeta )
 -\tau (z ))^2}  -\f{f(z )}{2 \pi i }  \int_{\pd D} \f{ d\ov{\tau(\zeta )}  }{\tau (\zeta )
 -\tau (z )}.
\eeq
On $\pd D$, write $d\ov{\tau(\zeta )}=a_0(\zeta )\, d {\tau(\zeta )}$. By (ii),
 we know that
the last integral  in \re{gztf} is in $\B^{\all,0}(\ov\D_r^+,P)$.
Name the first integral in \re{gztf} by $\tilde g(z )/{(2\pi i)}$. That  $\tilde g$ extends continuously follows
from the continuity of $f$ and $|f(\zeta)-f(z)|/{|\tau(\zeta)-\tau(z)|^2}\leq C|\zeta-z|^{\all-2}$. Write
\aln
\tilde g(z_2 )-\tilde g(z_1 )&= \int_{D} \f{(f(z_1 )-f(z_2 ))\,
d\tau(\zeta )\wedge d\ov{\tau(\zeta ) }}{
(\tau (\zeta )-\tau (z_2 ))(\tau (\zeta ) -\tau (z_1 ))}\\
&\quad + \int_{D}
   \f{(f(\zeta )-f(z_2 ))(\tau(z_2 )-\ta(z_1 ))}{(\tau (\zeta )
 -\tau (z_2 ))^2(\ta(\zeta )-\tau(z_1 ))}\, d\tau(\zeta )\wedge d\ov{\tau(\zeta ) }\\
 &\quad  + \int_{D}
   \f{(f(\zeta )-f(z_1 ))(\tau(z_2 )-\ta(z_1 ))}{(\tau (\zeta )
 -\tau (z_1 ))^2(\ta(\zeta )-\tau(z_2 ))}\, d\tau(\zeta )\wedge d\ov{\tau(\zeta ) }.
\end{align*}
%
%
%
The last two integrals can be estimated by a standard
argument for H\"older estimates,  bounded in absolute
value by
$ C_\all\|f\|_{\all,0}|z_2-z_1|^\all.$ The first integral can be rewritten as
the product of $
 f(z_1 )-f(z_2 )
$ and $\cL I$ for
\aln
\cL I
& =\f{1}{\tau (z_2 )-\tau (z_1 )}\int_{D}\biggl\{ \f{d\tau (\zeta )\wedge d\ov{\tau (\zeta )}}
{\tau (\zeta )-\tau (z_2 )}- \f{d\tau (\zeta )\wedge d\ov{\tau (\zeta )}}
{\tau (\zeta )-\tau (z_1 )} \biggr\}
\\
 &=2\pi i\f{\ov{\tau (z_2 )-\tau (z_1 )}}{\tau (z_1 )-\tau (z_2 )}
+
\f{1}{\tau (z_1 )-\tau (z_2 )}\int_{\pd D}\biggl\{ \f{ \ov{\tau (\zeta )}\,
d\tau (\zeta )}
{\tau (\zeta )-\tau (z_2 )}- \f{ \ov{\tau (\zeta )}\,
d\tau (\zeta )}
{\tau (\zeta )-\tau (z_1 )}\biggr\}.
\end{align*}
A  derivative of  $\int_{\pd D}\f{ \ov{\tau (\zeta )}\,
d\tau (\zeta )}
{\tau (\zeta )-\tau (z )}$ is $\pd_z\tau(z) \int_{\pd D}\f{
d \ov{\tau (\zeta )}}
{\tau (\zeta )-\tau (z )}$ which,  by  (ii), is bounded. By the mean-value-theorem,
the last term in $\mathcal I$ is bounded.
This shows that $S_0f\in\cL C^{\all,0}(\ov\D_r^+,P)$.

For higher order derivatives, we  transport derivatives to $f$.
Define $ h_*(\tau(z ) )=h(z )$
and    $\om(t)=\tau(\cdot,t)(D)$.
Let $C_*=C_{\pd\om(t)}$, $T_*=T_{\om(t)}$ and  $S_*=S_{\om(t)}$.
 Rewrite \re{gztf} as
$
  g_*(\tau )=S_{*}f_*.
$
Integrating by parts, we obtain
\aln
 g_*(\tau ) & = \f{1}{2\pi i}\int_{\om(t)}
 \f{\pd_{\fta}f_*(\fta )}{\fta-\tau}
  \, d\varsigma\wedge d\ov\varsigma
  -  \f{1}{2\pi i}\int_{\pd\om(t)}\f{f_*(\fta )}{\fta-\ta}\, d\ov\fta.
\end{align*}
On $\pd\om(t)$, we write $d\ov\ta=a(\ta,t)\, d\ta$ with $a
\in\cL C^{k+\all,0}(\pd D,P)$.
Taking derivatives, we get
$$
\pd_{\ov\ta}S_*f= \pd_\ta f_*,\quad
 \pd_\ta S_*f_*=S_*\pd_\ta f_*-\pd_\ta C_*af_*.
$$
Using the last formula $k$ times, we get
\aln
 (\pd_{\tau})^{k}S_*f_*  & =S_* \pd_\ta^kf_*-
 \sum_{0\leq j<k} \pd_\ta^{k-j}  C_*a\pd_{\ta}^jf_*.
\end{align*}
We return to the $z$ coordinates. Let $\tilde a(z,t)=a(\ta(z,t),t)$.
 Let $\pd_z^K$ be  a derivative
in $z,\ov z$ of order $k$.
 Let $\pd^{(j)}$ denote derivatives   of orders $\leq j$.
 Then
\al\label{s0ch}
 \pd_z^K S_0 f(z)&=
p^1 \bigl(\pd_z^{(k)}\tau\bigr) \cdot\bigl(\pd_\ta^{(k)}S_*f_*,\pd_\ta^{(k)}f_*\bigr)\circ\ta\\
&  =\sum_{0\leq j\leq k}
p^2_{ j }\bigl(\pd_z^{(k)}\tau\bigr) \cdot\bigl(S_*\pd_\ta^{(k)}f_*,
\pd_\ta^{(k-j)}C_*(a \pd_\ta^{(j)}f_*),\pd_\ta^{(k)}f_*\bigr)\circ\ta
\nonumber\\
&   =\sum_{0\leq j\leq k}
q^1_{j, k }   \cdot\bigl(S_0q_{k}^2
\pd_z^{(k)}f,
\pd_\ta^{(k-j)}C_0(\tilde a q_{j}^3 \pd_\ta^{(j)}f ),\pd_\ta^{(k)}f\bigr).
\nonumber\end{align}
Here integral operator  $S_0 $ is over the domain $D$. And
 $ C_0$ is over $\pd D$.  $p_j ^i $ are vectors of
 polynomials, and $q_{l,j}^1,
 q_j^i$ are matrices of polynomials in $(\det\tau')^{-1},\pd_z^{(j)}\tau$.

That $S_0f\in\B^{k+\all,0}(\ov\D^+,P)$ follows from the assertion for $k=0$ and (ii).

   Note that $\pd_\ta T_*=S_*$ and $\pd_{\ov\ta} T_*=\I$.
Thus, $T_0f\in\B^{k+1+\all,0}(\ov\D_r^+,P)$ by (ii),    the product rule,
and the  chain rule as used
in \re{s0ch}. \end{proof}

\setcounter{thm}{0}\setcounter{equation}{0}
\section{Proof of the higher dimensional result
}\label{sec4}

Let $\Delta_r^{n}, \Del_r^{2n-1}, \Del_r^{2n}$ be the polydiscs of radius
$r$ in the $x$-subspace,   hyperplane $y_n=0$,
and   $\rr^{2n}$, respectively.

In this section, we will prove \rt{combvn}.  In view of \rl{jjv} and Example~\ref{ex+},
it is worth stating a more general result. This will also make the proof transparent.


\th{combvn2}
 Let $k\geq   4$ be an integer.  Let $\Om_1,\Om_2, M,\alpha$ be as in \rt{combvn} with
$M\in\cL C^{k+1+\all}$.  For $i=1,2$, let $
J^i
$
be an almost complex structure of class
$\cL C^{k+\all}(\Om_i\cup M)$
on $\Om_i\cup M$.
    Suppose that at each point $p\in M$ there is a tangent
vector $v_p\in T_pM$ such that $J^1_pv_p, J^2_pv_p$ are in the same connected component
of $T_p\rtn\setminus T_pM$.
Let $f\in\cL C^0(\Om_1\cup\Om_2)$ be a continuous  function  on $\Om_1\cup M\cup \Om_2$
such that
  $ (\pd_{x_j}+\sqrt{-1}J^i\pd_{x_j})f$ and $ (\pd_{y_j}+\sqrt{-1}J^i\pd_{y_j})f$, defined on $\Om_i$,  extend to
 functions in $ \cL C^{  k }( \Om_i\cup M)$  for $i=1,2$, $j=1,\ldots, n$.
 Then   $f$ is of class $\cL C^{k-3+\beta}(\Om_1\cup M)$
 for all $\beta<1$.\eth

 Notice that no integrability
condition
 is assumed.  A by-product of our proof is  $f\in
 \cL C_{loc}^{k-3+\beta}(\Om_1)$ for all $\beta<1$
 when  $k\geq3$. (Of course the assumptions on $f, J^2$ for $M,\Om_2$
 are not needed.)

\smallskip

The main ingredients of the proof are in the following.

Step 1.  We will show that
the Fourier transform
of $f$ on lines $L$ in $M$ decays
  in the $\xi$-variable. To use the differential equations for $f$,  lines $L$ need to be
   transversal to
  the complex tangent vectors of $M$ of both   structures.
  Two
almost complex structures yield decay of the Fourier transform at opposite rays.
This is the only place we need  both   structures.
This gives us  smoothness of $f$ on $M$.

Step 2. We will   obtain smoothness
of $f$ on each side of $M$ (up to the boundary) via the one-sided almost complex
structure.
We attach a family of holomorphic  discs to $M$ with respect to the  structure.
Such a disc will have regularity as good as the structure provides.
This is achieved by extending the   structure
to  a neighborhood of $M$.
The regularity of $f$ on $M$ yields uniform bounds of pointwise derivatives
of $f$ along the  discs up to their boundaries in $M$.

Step 3. Let $\Om^+=\Om_1$. After   obtaining   smoothness of $f$ on families of
  discs in $\Om^+\cup M$, we conclude
the smoothness of $f$ on   of $\Om^+\cup M$ via \rp{ftva}.

\medskip

We now carry out   details. We need a preparation for Step 1.

\smallskip

\noindent{\bf Step 0. Match approximate $J$-holomorphic half-discs   in $M$.}
%

We may assume that $M$ is $\Del^{2n-1}\times0$, $\Om^+=\Del^{2n}\cap\{y_n>0\}$
and $\Om^-=\Del^{2n}\cap\{y_n<0\}$.   Let $0<r<1$ be sufficiently small.
   By the assumption, there is a  vector $v_0\in  T_0M$ such that
the vectors $J^1_0v_0,  J^2_0v_0$ are
transversal to $T_0M$ and are in $\Om^+$. Thus
the line segments $tJ^1_0v_0,  tJ^2_0v_0$ ($0<t\leq1$) are
transversal to $M$ and are in $\Om^+$, by shrinking $v$ if necessary.
Here we have identified $\rtn$ with $T_p\rtn$  by sending $v$ to the tangent
vector of $p+tv$; consequently, $J^i_p$   acts on $\rtn$ linearly.
  Let $\e>0$ be sufficiently small, let  $p\in M, v\in T_pM$ satisfy $|p|<\e$ and  $|v-v_0|<\e$.
By transversality,   $p+tJ_p^1  v$ and $ p+ tJ^2_p  v$ are in
  $\Om^+$ for $0<t\leq1$. 
    Define
  $$
 L= L(v,p)=\{p+sv\colon -2<s<2\}\subset M.
  $$
  Let $e_1,\ldots, e_{2n-1}$ be the standard basis of $\rr^{2n-1}$.
We find an affine coordinate map $\phi$ on   $\rtn$ such that $\phi(p)=0$, $\phi(p+v)=e_1$,
and  $\phi(p+v_j)=  e_j$.
We may also assume that the norms of $\phi$ and $\phi^{-1}$ have an upper bound independent
of $p$, $v$. In what follows, all constants are independent of $p,v$.
\rp{ftva} (ii) will be used for this family of $\phi$ (with $p=0$) depending on parameter $v$ with $v_0$
to be chosen.

 We want to apply \rl{flat} to $L(v,p)$. Here $v,p$ are parameters and we suppress them in all
 expressions.
 For the above $L(p,v)$, we attach an approximate $J$-holomorphic curve $u^1$   of class $\cL C^{k+1+\all}$
 such that
\ga\label{xjyj}
du^1(\pdoz)=D^1(z)\cdot X^1(u^1(z))+F^1(z)\ov{X^1(u^1(z))},\\
\nonumber 
|F^1(z)|\leq C|y|^{k+\all}, \quad (x,y)\in Q\df (-1,1)\times (0,\e).
\end{gather}
We have an analogous $u^2$ on $\Om^-\cup M$. We have
$$
u^1(x,0)= p+xv=u^2(x,0)  \  \text{on $[-1,1]$}.
$$
We know that $u(x,0)$ is contained in $M\subset\ov\Om_1^+\cap\ov\Om_1^-$
for $|x|<1$. When $p=0$ and $v=v_0$, we have $du^1(0)(\pd_x)=v_0$ and
  $du^1(0)(\pd_y)=  J_0^1du^1(0)(\pd_x)=J_0^1v_0$ is contained in $\Om_1^+$,
  $-  J_0^2v_0$ is contained in $\Om_1^-$
and both are transversal to $ M$. Thus,
  \ga\label{distm-}
   u^1(x,y)\in\Om^+, \  (x,y)\in Q; \quad
  u^2(x,y)\in\Om^-, \  (x,y)\in -Q.
  \end{gather}
The above hold for $v=v_0$ and $p=0$.
Since the derivatives of $u$ are continuous in $p,v$, the above hold for $|p|<\e$ and $|v-v'|<\e$. And for
a constant $C>1$ independent of $p,v$,
\ga 
\dist(u^i(x,y),M)\geq |y|/C,\quad (x,y)\in (-1)^{i-1}Q.
\label{distm}\end{gather}

\noindent{\bf Step 1. Uniform bound of Fourier transform of $f$ on
 transversal lines $L$ in $M$.}

\nopagebreak
In this step and the next, we will assume that $f$ is $\cL C^1$ on $\Om_1\cup\Om_2$. We will
verify this interior regularity
in the final step.

Fix $k$. Recall from Step $0$ that $M$ is contained in $\rr^{2n-1}$. Let $v_0,\e$
be as in Step~$0$.
By \rl{det}, there exist $d$  vectors $v_j$ in $\rr^{2n-1}$ such that
\eq{pdxy}
(\pd_x,\pd_{y'})^\all=\sum_{1\leq j\leq d}
 c_{\all, j}(v_j\cdot(\pd_x,\pd_{y'}))^{|\all|}, \quad 1\leq|\all|\leq k.
\eeq
Here
  $|v_j-v_0|<\e$.
Recall the line segment
 $L=\{p+sv_i\colon -1\leq s\leq 1\}$ with $p\in M, |p|<\e$.
 Fix such an $L$ and denote its tangent vector $v_i$ by $v$.

\medskip

Note that when $\e$ is sufficiently small, $L$ has length $>|v_0|/2$.
Let $\chi_0$ be a cutoff function on $M$ with compact support in $\Del_{|v_0|/{(4n)}}^{2n}\cap M$.
Then $\chi_0|_L$ has compact support.
We will show that the Fourier transform of $\chi_0 f$ on $L$ satisfies
\eq{foli}
(1+|\xi|)^{k-1+\alpha-\beta}\left|\widehat{ \chi_0 f |_L}(\xi)\right|<C_\beta\eeq
 for all $\beta>0$. Here $C_\beta$
 will be independent of $p$,    $v_1,\ldots, v_d$.
We will verify \re{foli} for $\xi=|\xi|v$, using $X_j^1f=g_j^1$ on $\Om^1$
  with $g_j^1\in\cL C^{k}(\Om^+\cup M)$.
For $\xi=-|\xi|v$, we   use $X^2_jf=g_j^2$ on $\Om^2$ with $g_j^2\in\cL C^{k}(\Om^-\cup M)$.

We now use approximate $J$ holomorphic curves $u^1,u^2$ defined in Step~$0$.
  We   drop the superscript in $u^1,g_j^1, a^1_{jk}$, etc.

Applying Whitney's extension   (\rl{whit}), we
extend $\chi_0\circ u(x,0)$ to    $\chi\in\cL C^{\infty}(Q)$ which has
compact support in each $(-1,1)\times\{y\}$. Moreover,
$|\pdoz\chi(x,y)|\leq C|y|^{k+\all}$.
For brevity, denote $f\circ u,  g_j\circ u, (\pd_{z_j}f)\circ u$
by $f, g_j, h_j$.  Combining with \re{xjyj},
we get  on $Q$
\ga\label{eqfchi}
\pd_yf(x,y)=i\pd_xf(x,y)-2iD(x,y)\cdot g(u(x,y))-2iF(x,y)\cdot \ov{X(u)}  f,\\
\pd_y\chi(x,y)=i\pd_x\chi(x,y)+E(x,y).\nonumber
\end{gather}
Here
$(|E|+|F|)(x,y) \leq C|y|^{k+\all}$. And $D, E,F$ are in $\cL C^{k+\all}(\ov{Q})$,
and  $g$ is in $\cL C^{k}(\ov{Q})$.

\medskip

In what follows, as required by \re{foli}
constants do not depend on $L$, $p, v_j$.

By \re{distm-},  $u(x,y)$ is in $\Om_1^+$ for  $|x|<1$,   $0<y<\e$.
Define
$$
\lambda(\xi,y)=   \int_{\rr} (\chi f)(x,y)e^{-i(x-iy)\xi }\, dx, \quad y\geq0.
$$
Note that
 $\widehat{\chi f|_L}(\xi)\equiv\lambda(\xi,0)=
\lambda(\xi,\eta)-\int_0^\eta\pd_y\lambda
(\xi,y)\, dy$. By \re{eqfchi}, we obtain
\aln
 \pd_y \lambda(\xi,y)&=\int_{\rr} \bigl(i \pd_{x}(\chi f)(x,y )-(\chi f)(x,y )
\xi\bigr)e^{-i(x-iy ) \xi}\,
dx\\ &\quad -
2i\int_{\rr} (g(u)\cdot D  \chi)(x,y)e^{-i(x-iy ) \xi}\,
dx
+\int_{\rr}( f(u)E)(x,y)e^{-i(x-iy)\xi}\, dx\\
&\quad -
2i\int_{\rr}   \chi F(x,y)\cdot(\ov Xf)(u(x,y))e^{-i(x-iy)\xi}\, dx.
\end{align*}
By   integration by parts,   the first integral is zero. Since $g(u(x,y)), D(x,y)\in\cL C^{  k}$
and $\eta\xi\geq0$, the second,
via integrating by parts $k$ times, is less than $C(1+|\xi|)^{-k}$.
The third  is bounded by $C|E(x,y))|\leq Cy^{k+\all}$. We now estimate
the last integral. This amounts to controlling the blow-up of  derivatives of $f$
at $u(x,y)$.  By \re{distm}, $ \Om^+$ contains
  $\Del_{y/C}^{2n}(u(x,y))$. To apply \rp{corZ} to the latter, we need a domain of fixed size.
    Let  $\psi(\zeta)=u(x,y)
  +\zeta y/C$. So   $\psi^{-1}$ transforms $J, X_j$ into $\hat J,
  \hat X_j=C^{-1}yd\psi^{-1}X_j$. On $\Del^{2n}$,
  we have
  $$
  \hat X_j=\sum_{1\leq k\leq n}(b_{jk}\circ\psi\pd_{z_j}+a_{jk}\circ\psi\pd_{\ov z_j}).
  $$
Let $ A'=(a_{jk}\circ\psi),   B'=(b_{jk}\circ\psi)$.  It is easy to see that on $\Del^{2n}$,
 $\inf {\tiny \begin{vmatrix}  {  B'} \ {  A'}\\
  { \ov {  A}'}\ { \ov {  B}'}\end{vmatrix}} \geq 1/C$ and $|(  A',  B')|_{k+1+\all}\leq C$
  for some constant independent of $v,p$.
Applying \rp{corZ} to $\{\hat X_j\}$, we get
 a $J$-holomorphic curve
$\hat u\colon\D_r\to\Del^{2n}$ with $\hat u(0)=0$,  $d\hat u(0)\pd_\zeta=\pd_{x_m}-i\hat J_{0}\pd_{x_m}$.  Here  $r>0$ is a constant independent of $y$.
 Then $\tilde u(\zeta)=\psi \circ \hat u(C\zeta/y)$ is
     $J$-holomorphic in $ J$.  We have $\tilde u\colon \Delta_{  y/c}\to \Om_1^+$,
$\tilde u(0)=u(0)$, and
 $$d\tilde u(0)(\pd_\zeta)=\pd_{x_m}-iJ_{u(0)} \pd_{x_m}. 
 $$
 So $d\tilde u(0) \pd_{\ov\zeta}=\pd_{x_m}+i  J_{u(0)} \pd_{x_m}.$
A direct computation shows that the first and second order derivatives of $\tilde u$
are bounded by $C$,
   $C/y$, respectively.
Since   $X_jf=g_j$ and $d\tilde u(\pd_{\ov\zeta})=\tilde D\cdot X(\tilde u)$,
then
$$ 
\pd_{\ov \zeta}(f(\tilde u(\zeta)))= g(\tilde u(\zeta))\cdot\tilde D(\zeta).
$$ 
Note that the derivative of $\tilde D$ is bounded by $C/y$.
By the Cauchy-Green identity, we have
$$
f(\tilde u(z))=\f{1}{2\pi i}\int_{  |\zeta|={ y/c}}\f{f(\tilde u(\zeta))}{\zeta-z}\, d\zeta
+\f{1}{\pi}\int_{|\zeta|<  y/c}
\f{ g(\tilde u(\zeta))\cdot\tilde D(\zeta)}{  z-\zeta}\, d\xi d\eta.
$$
  At $z=0$,   derivatives of the first integral are bounded
by $C/y$. Write
 $g(\tilde u(\zeta))\cdot\tilde D(\zeta)$ as $h_1(\zeta)+h_2(\zeta)$. Here
  $\cL C^1$ norms of $h_1,h_2$ are bounded by $C/y$, $h_1(\zeta)=0$
on $|\zeta|<y/{(4c)}$,  and $h_0(\zeta)=0$ on $|\zeta|>y/{(2c)}$. The derivatives of
the integral involving $h_1$ are bounded by $C$ at $z=0$. After applying
   translation $\zeta'=\zeta-z$,
  the integral  involving $h_0$  has bounded derivatives
  at $z=0$ too.
We obtain  $|\pd_{x^m}f(u(0))|=|\pdoz f(u(0))+\pd_zf(u(0))|\leq C/{y}$.
   Thus,
$$
|\pd_y \lambda(\xi,y)|\leq
  C\left(y^{k-1+\all}
+ (1+|\xi|)^{-k}\right).
$$
We also have $|\lambda(\xi,\eta)|\leq Ce^{-  \eta\xi }\leq C_L|\eta\xi|^{-L}$. We may assume that $\xi\geq1$.
Choose   $\eta=1/{(C|\xi|^\all)}$. Then $\la(\xi,0)=\la(\xi,\eta)-\int_0^\eta\pd_y\la(\xi,y)\, dy$ satisfies
\eq{lade}
  |\lambda(\xi,0)|\leq C (1+|\xi|)^{-(k+\all)}
\eeq
for $ \xi\geq0$. Reason by using $  X_j^2f=  g^2_j$ for $y_n\leq0$, $u^2$,
 and by replacing $e_n$ with $-e_n\in\rr^n$. We get \re{lade} for $\xi\leq0$ and hence
  for $-\infty<\xi<\infty$.

By the Fourier inversion formula,
\gan
\chi f(p+xv)= \f{1}{2\pi}\int_\rr \la(\xi,0)e^{i\xi x}\, d\xi,\\
\pd_x^{k-1}(\chi f(p+xv))= \f{1}{2\pi}\int \la(\xi,0)(i\xi)^{k-1}e^{i\xi x}\, d\xi.
\end{gather*}
This shows that $\chi f(p+xv)\in\cL C^{k-1}$. By the mean-value-theorem and \re{lade},
$$
|\pd_{x}^{k-1}(\chi f)(p+x_2v)-\pd_{x}^{k-1}(\chi f)(p+x_1v)|
\leq C\int\f{|x_2-x_1|^{\all'}|\xi|^{\all'}}{(1+|\xi|)^{ 1+\all }}\,d\xi.
$$
For any $\all'<\all$, we have
   $
 | (\chi f(p+\cdot v))|_{k-1+\all'}<C_{\all'}.
 $
Therefore,
  $$
 \left|  \chi_0 f|_L\right|_{k-1+\all'}<C_{\all'},
 $$
where $L$ is any line which is
 tangent to one of $v_1,\ldots, v_d$ and passes through $p$ for $p\in M$ near the origin.
 For the $L$, we can find an affine diffeomorphism
 $R$ with $R(0)=0\in M$, sending
 $\Del^{2n-1}$ into $M$, such that
  $R(\cdot, t)$ are lines parallel to $L$ for $t\in\Del^{2n-2}$.
By \rp{ftva} (ii) and hypothesis \re{pdxy}, we get $\pd^{k-2}(  \chi_0 f)\in Lip\,(M)$.

\medskip

\noindent{\bf Step 2. Uniform bound of derivatives of $f$ on transversal $J$-holomorphic
curves.}

\nopagebreak

Fix $k$. By \rl{det}, there exist $N$ vectors $v_j\in\rr^{2n}$  with
  $|v_{j}|<1$ such that
\eq{pdxy2n}
(\pd_x,\pd_{y})^\all=\sum_{1\leq j\leq N}
 c_{\all, j}\bigl(v_j\cdot(\pd_x,\pd_y)\bigr)^{|\all|},\quad 1\leq |\all|\leq k.
\eeq
 By perturbing $v_j$, we may  assume that
$J_p(v_j\cdot(\pd_x,\pd_y))$ are not tangent
to $M$ at $p=0$ and hence in a neighborhood of $0\in M$.

\smallskip

We are given vector fields  $\{X_j\}$ defined on $\Om^+\cup M$. Recall that
$M$ is contained in $y_n=0$ and   $\Om^+$ is contained in $y_n>0$. Applying
 Whitney's theorem via restriction and then extension,
we may assume that $\{X_j\}$ is an almost
complex structure  defined in a neighborhood of $M$.
We now apply \rp{corZ} with the parameter set $P=\ov\D_{r_0}^{n-1}$.
There are diffeomorphisms $u_j, R_j$
of class $\cL C^{k+1+\all, k-1}$,
which maps $\D_{r_0}\times\D_{r_0}^{n-1}$ into
$\Om$ with
\eq{dujvj}
du_j(0,t)(\pd_{\xi})=v_j.
\eeq
Moreover,  $D_{j,r}(t)=u_j(\D_r,t)=R_j(\om_{j,r}(t),t)$ satisfy
$$
\D_{r/c_2}\subset \om_{j,r}(t)\subset \D_{c_2r}.$$
Also $R_j(0)=0$, $u_j(0)=0$, and
$$
du_j(\cdot,t)(\pd_{\ov\zeta})=J_0 (v_j\cdot(\pd_x,\pd_y))-i(v_j\cdot(\pd_x,\pd_y)).
$$
We choose $r<r_0$ sufficiently
small so that   various compositions in $u_j,R_j$ are well-defined.
Also $\om_{j,r}^+(t)=\om_{j,r}(t)\cap\{y>0\}$ satisfies $R_j(\om_{j,r}^+(t),t) \subset
D_{j,r}^+(t)=D_{j,r}(t)\cap\Om^+$. Write $(\tilde D_{j,r}^+(t),t)=u_j^{-1}(D_{j,r}^+(t))$.
We apply the Cauchy-Green formula for $f\circ u_j(\cdot,t)$. Then
$$
f(u_j(\tilde z,t))=\f{1}{2\pi i}\int_{\pd\tilde D_{j,r}^+(t)}
\f{f(u_j(\tilde \zeta,t))}{\tilde \zeta-\tilde
z}\, d\tilde\zeta+\f{1}{2\pi i}\int_{\tilde D_{j,r}^+(t)}\f{\pd_{\ov{\tilde\zeta}}
f(u_j(\tilde \zeta,t))}{\tilde \zeta-\tilde
z}\, d\ov{\tilde\zeta}\wedge d\tilde\zeta.
$$
Set $(\tilde z,t)=u_j^{-1}\circ R_j(z,t)
=(\tau_j(z,t),t)$. By $du_j(\pd_{\ov\zeta})=D(\zeta)\cdot X(u_j(\zeta))$ and $X_jf=g_j$,
we get  
$\pd_{\ov\zeta}f(u_j)=D\cdot g $ and
\aln
f(R_j(  z,t))&=\f{1}{2\pi i}\int_{\zeta\in
\pd\om_{j,r}^+(t)}\f{f(R_j(  \zeta,t))}{\tau_j(\zeta,t)-\tau_j(z,t)}\, d\tau_j(\zeta,t)
\\ &\quad+\f{1}{2\pi i}\int_{\om_{j,r}^+(t)}\f{g(R_j(\zeta,t))
\cdot D(\tau_j(\zeta,t),t)}{\tau_j(\zeta,t)-\tau_j(z,t)}
\, d\ov{\tau_j(\zeta,t)}\wedge d\tau_j(\zeta,t).
\end{align*}
Recall that $R_j(\cdot,t)$ send $[-r/c,r/c]$ into $\pd D_j^+(t)\cap M$. Since $u\in\cL C^{k+1+\all,k-1}$,
it is easy to see that $D(\tau_j(\zeta,t),t)$ are in $\cL C^{k+\all,0}$.
Applying \rl{ctva}, we get $|f(R_j)|_{k-2+\beta,0}<C_\beta$ on $D_{r/c_*}$ for any $\beta<1$.

\medskip

\noindent{\bf Step 3.  Smoothness of $f$ via families of $J$-holomorphic curves.}
\nopagebreak

By the end of Step 2, we know that the $\cL C^{k-2+\beta}$ norms
$f\circ R_j(\cdot,t)$ on  $\ov\D_{r}^+$ are  bounded
in $t\in\D_{r}^{n-1}$ when $r$ is small, and that $\pd_{t_1}^{k-2}(f\circ R_j(t))$
is continuous. Here $R_j(0)=0\in M$.
Applying the Whitney extension, we extend $f(z)$ to $\Om^-$ such that
it has class $\cL C^{k-2+\beta}$ on $\Om^-\cup M$.
Then  $f\circ R_j$
is of class $\cL C^{k-2+\beta}$ on $\D_r$. Therefore,
for the extended function $f$, the  $\cL C^{k-2+\beta}$ norms of $f(R_j(\cdot, t))$ on $\D_r$
are  bounded and $\pd_{\xi}^{k-2}(f\circ R_j(\xi+i\eta, t))$ is continuous.
We want to apply \rp{ftva} (iii) to a family of diffeomorphisms $\tilde
R_j\in\cL C^{k-1}$.
This is achieved easily by taking $\tilde R_j(t_1,t_2,t')=u_{j}(t_1+it_2, t' )$
and treating $(t_2,t')\in\rr^{2n-1}$ as parameters.
  By \re{pdxy2n}-\re{dujvj}, we conclude $f\in \cL C^{k-3+\beta}(\Om^+\cup M)$ for all $\beta<1$.

\smallskip
To finish the proof, we need to remove the assumption stated at
the beginning of Step 1 that $f$ is $\cL C^1$ on $\Om_i$.
We also  address the comment made after
\rt{combvn2}
on the interior regularity of $f$.
 Let $J=J^i\in\cL C^{k+\all}$  and $\Om=\Om_i$.
Here we only need $k\geq2$.  Let $v_j$ satisfy
\re{pdxy2n}. By \rp{corZ} we find a $\cL C^{k+1+\all,k-1}$ diffeomorphism $u_j(\zeta,t)$ defined in neighborhood of $0\in\Om$
such that $\zeta\to u_j(\zeta, t)$  is  $J$-holomorphic for fixed $t\in\D_\e^{n-1}$
 and $du_j(0)(\pd_\xi)=v_j$. Drop the subscript $j$ in $u_j$.
 Then $u^{-1}$ defines a $\cL C^{k-1}$ coordinate system for a neighborhood
   of  the origin. And $u^{-1}$
  transforms $J$ into $\hat J$. Now $\D_\e\times t$ are $J$-holomorphic discs
 in $\hat J$  for $|t|<\e$
 and $\e$ small. Thus we can take $\hat X_1=a(\zeta,t)\pd_{\ov \zeta}+b(\zeta,t)\pd_{\zeta}$ with $a,b\in
 \cL C^{k+\all, k-1}$  and $t$ being parameters. Now $\hat f=f\circ u$
 satisfies $\hat X_1\hat f=\hat g_1\in\cL C^{k,k-1}$. Here $\hat X_1f=\hat g_1$ holds in the sense of distributions.
 We want to show that  when restricted on
  $\D_\e\times t$, $\hat X_1\hat f=\hat g_1$ still holds as distributions.
  To verify it,  fix a test function $\phi$  on $\D_\e$ and take a sequence of test functions $\phi_i$ in $ \cc^{n-1}$
 such that $\int_{\cc^{n-1}} \phi_i=1$ and $\supp\phi_i
 \subset B_{1/j}(t)$. Note that
  the formal adjoint
 $\hat X_1^*$ does not contain derivatives in $t$-variables and satisfies
 $$
 \int_{\cc^n} \hat g_1\phi\phi_i=\int_{\cc^n}\hat  f\hat X_1^*(\phi\phi_i)=\int_{\cc^n}\hat  f\phi_i\hat X_1^*(\phi).
 $$
Since all functions in the integrands are continuous, letting $i$ tend to $\infty$ yields
$$
\int_{\cc}\hat g_1(\cdot, t)\phi=\int_{\cc}\hat  f(\cdot,t) \hat X_1^*(\phi).
$$
We have proved that $\hat X_1\hat f=\hat g_1\in\cL C^{k,k-1}$ in the sense
of distributions and the coefficients of $\hat X_1$ are in $\cL C^{k+\all,k-1}$.
Reasoning as at the end of Step 2
 by the Cauchy-Green identity, we see that   $f\circ u_j \in\cL C^{k+\all,0}$.  Now $f\circ u_j\in\cL C^{k-2,0}$,
$u_j\in \cL C^{k-1}$ and
 \rp{ftva} (iii) implies that $f\in\cL C^{k-3+\beta}$ for all $\beta<1$.
 By $k\geq4$, the   proof of \rt{combvn2} is complete.

\setcounter{thm}{0}\setcounter{equation}{0}
\section{One-dimensional results
}\label{sec5}
Throughout this section,
   $\Om$ is a bounded open set in  $\cc$, and
$P$ is the closure of a  bounded open set in a  euclidean space. We assume that
      two points $a,b$ in $\ov\Om\times  P$ can be connected
by a smooth curve in $\ov\Om\times   P$ of length at most $C|b-a|$.

We start with  the existence of isothermal coordinates with
parameter.   Recall that  a diffeomorphism $\var$ is said to transform a vector field $X$ into
$\tilde X$ if locally $d\var(  X)=\mu \tilde X$. Denote by $\B_{loc}^{k+ \all,j}
(\Om\cup\gaa,P)$ the set of functions which are in $\B^{k+ \all,j}
(K,P)$ for any compact subset $K$ of $\Om\cup\gaa$.
\pr{isop} Let $\Om$ be a domain in $\cc$ and
$P$ be an open set in a  euclidean space.
Let $a\in\B^{k+\all,j}(\Om,P)$ satisfy $|a|_{0,0}<1$.  If $x\in\Om$ there exist a neighborhood $U$
of $x$ and a map $\var\in\B^{k+1+\all,j}(U,P)$ such that $\var(\cdot,t)$ are diffeomorphisms
which map $U$ onto their images and  $\pdoz+a(z,t)\pd_z$ into  $\pdoz$.
\epr
\begin{proof} Fix $x=0\in\Om$. Let $\var(z,t)=z-a(0,t)\ov z$. Then $z\to\var(z,t)$ is invertible
and transforms $\pdoz+a(z,t)\pd_z$ into $\pdoz+\tilde a(z,t)\pd_z$ with $\tilde a(0,t)=0$. We still have
$|\tilde a|_{\all,0}<\infty$.
 Let $\chi$ be a smooth function on $\D$ which has compact support and equals $1$ on $\D_{1/2}$.
 Applying a dilation and replacing $\tilde a$ by $\chi\tilde a=b$ we achieve $|b|_{\all,0}<\e_\all$ on $\D$
for $\e_\all$ in \rl{isa}. Set
$
f=-(\I+T  b\pd_z)^{-1}Tb.
$
On $\D$ we have $f\in\B^{k+1+\all,j}$ and $|f|_{1,0}\leq C_\all|(\I+T  b\pd_z)^{-1}|_{1+\all,0}|b|_{\all,0}$.
With the dilation for $\tilde a$, $|f|_{1,0}$ can be arbitrarily small. Therefore
$z\to z+f(z,t)$ are indeed diffeomorphisms. Since
$z+f(z,t)$ is annihilated by $\pdoz+\tilde a\pd_z$, it transforms $\pdoz+\tilde a(z,t)\pd_z$
into  $\pdoz$.  \end{proof}

It is important that the above classical result (for the non-parameter
case) allows one
to interpret  $\pdoz f+a\pd_zf=g$
when $a$ is merely $\cL C^\all$. Let $w=\var(z)$ be a  local $\cL C^{1+\all}$
diffeomorphism  such that
$d\var(\pdoz+a\pd_z)=\mu(w)\pd_{\ov w}$.
Then $\pdoz f+a\pd_zf=g$ holds in the $w$-coordinates,
 if $\pd_{\ov w}(f\circ\var^{-1})= g\circ\var^{-1}(w)/{\mu(w)}$ holds in the sense of
distributions.
We restate \rp{combv} in a parameter version.
\pr{combv+}   Let $0<\all<1$ and $k\geq j\geq0$ be   integers.
Let $\gaa$ be an embedded curve in $\cc$ of class
$\cL C^{k+1+\all}$. Let $\Om_1,\Om_2$ be disjoint open subsets of $\cc$ such that
  both  $\pd\Om_1,\pd\Om_2$ contain
$\gaa$ as relatively open subsets.   Assume that
 $a_i \in\B^{k+\all,j}(\Om_i\cup\gaa,P)$
 satisfies $
|a_i|_{0,0}<1$ on $(\Om_i\cup\gaa)\times P$.
Let $f\in\B^{0,j}( \Om_1\cup\gaa\cup\Om_2,P), b_i\in\B^{k+\all,j}(\Om_i\cup\gaa,P)$   satisfy
$$
\pdoz f+a_i \pd_zf=b_i\ \text{on $\Om_i$},\qquad i=1,2.
$$
Then   $f\in\B_{loc}^{k+1+\all,j}
(\Om_i\cup\gaa,P)$.
\epr
\begin{proof}
As in the proof of \rt{combvn} in section~\ref{sec4}, we may assume that $\gaa$ is the $x$-axis
and $\Om_1=\D_r^+, \Om_2=\D_r^-$. In the following, all functions $a_i,b_i$, etc.~are
defined on $\Om_i$ for some $r>0$ and we will take smaller values for $r$ for a few times.
Apply the Whitney extension theorem with parameter (\rl{whit}).
We first find a function $\phi_i\in\B ^{k+1+\all,j}$ with $\phi_i(\cdot,t)\in\cL C^2(\Om_i\cup\gaa)$ such
that $\phi_i(x,0,t)=x$ and $\pdoz\phi_i+a_i\pd_z\phi_i=O(|y|^{k+\all})$.
 Then $\phi_i$ sends
$\pdoz+a_i\pd_z$ into $\mu_i(\pdoz+\tilde a_i\pd_z)$.  Replace $f, a_i$ by $f\circ\phi_i^{-1},\tilde a_j$
on $\ov\Om_i$.
Therefore, we may assume that $a_i(z,t)=O(|y|^{k+\all})$.
Define $a=a_i$ on $\ov\Om_i$. Then $a$ is of class $\B^{k+\all,j}(\Om_1\cup\gaa\cup\Om_2,P)$.
Set $X=\pdoz+a(z,t)\pd_z$.

Next, we find $g_i\in\B ^{k+1+\all,j}$ on $\ov\Om_i$ so that
$$
X g_i-b_i=O(|y|^{k+\all}),\quad g_i(x,0)=0.
$$
Replace $f$ by $f-g_i$ on $\ov\Om_i$.
Therefore, we may assume that $b_i(z,t)=O(|y|^{k+\all})$.
 Define $b=b_i$ on $\ov\Om_i$. Then $b$ is of class $\B ^{k+\all,j}(\Om_1\cup\gaa\cup\Om_2,P)$.

By \rp{isop}, there are
diffeomorphisms $\psi(\cdot,t)$ with $\psi\in\B^{k+1+\all,j}(\D_r,P)$, which   send  $X$ into   $\mu\pdoz$
with $\mu\in \B^{k+\all,j}(\D_r,P)$.
Then $\pdoz (f\circ\psi^{-1})=b\circ\psi^{-1}/\mu$. Let $h=T_{\D_r}(b\circ\psi^{-1}/\mu)$ where
$r$ is sufficiently small. Then $h\in\B^{k+1+\all,j}$.
Now $f\circ\psi^{-1}-h$ is holomorphic  away from $\psi(\gaa)$,  continuously up to the $\cL C^1$ curve
 $\psi(\gaa)$. Take a small disc $D_r$, independent of $t$ and
  centered at $p\in\psi(\cdot,t_0)(\gaa)$. By the Cauchy formula, we express
$f(\cdot,t)$   on $\Del_r$ via the Cauchy transform  on $\pd\Del_r$ when $t$ is in a small neighborhood
of $t_0$.
From $f\in\cL C^{0,j}$ and compactness  of $P$ we conclude
 $f\in\B^{k+1+\all,j}(\D_{r/2},P)$. Recall that $f$ is replaced by $f\circ\phi_i$.
The original  $f$ is in $\B_{loc}^{k+1+\all,j}
(\Om_i\cup\gaa,P)$.
\end{proof}

\le{teqs0} Let $D\subset\cc$ be a bounded domain with $\pd D\in\cL C^1$.
Suppose that $v\in
\cL C^1(D)$ and $b$ are continuous functions on $\ov D$.
Then    $v$ satisfies
\eq{ieq}
v +\T b  =0
\eeq
if and only if   it   satisfies
\ga\label{deq}
\pdoz v+b=0, \\
\label{bva}
\Ct v =\f{1}{2\pi i}\int_{\pd D}\f{v(\zeta)}{\zeta-z}\, d\zeta=0.
\end{gather}
Here three identities are on $D$.   Moreover,
\rea{bva} holds on $D$ if and only if
$v$ is the boundary value of a function that is holomorphic on $\cc\setminus\ov D$,
continuous on $\cc\setminus D$, and vanishing at $\infty$.
\end{lemma}
\begin{proof} Applying $\pdoz$ to \re{ieq} gives us \re{deq}.
On $D$,  $\Ct v=v-T\pdoz v$.
Applying $T\pdoz$ to   \re{ieq} and using \re{ieq} again, we get $v- T\pdoz v=0$.
Conversely, if $v$ satisfies \re{deq}, then
  $v+Tb=v-T\pd_{\ov z}v=\cL Cv$. The latter   is zero by \re{bva}.
Thus $v$ satisfies \re{ieq}.

  It is  a standard fact that when $D$ is a bounded domain with $\cL C^1$ boundary
and $v$ is continuous on $\pd D$, then  $\cL Cv(z-t n(z))-\cL Cv(z+tn(z))$
converges to $v(z)$ uniformly on $\pd D$ as $t\to0^+$. Here $n$ is the unit outer
normal vector of $\pd D$. Then \re{bva} implies that $\cL Cv$ is continuous
on $\cc\setminus D$ and agrees with $v$ on $\pd D$. That $\cL Cv$ vanishes
at $\infty$ is trivial. The converse follows from the Cauchy formula.
\end{proof}

Applying the above component-wise to the vector-valued functions, we get
\le{teqs1}  Let $D\subset\cc$ be a bounded domain with $\pd D\in\cL C^1$.
Suppose that $v\in
\cL C^1(D)$ and $b$ are vectors of $n$ continuous functions
on $\ov D$. \bppp
\item Let $A$ be an $n\times n$ matrix of continuous functions on $\ov D$.
Then    $v$ satisfies
\eq{ieq1}
v +\T (b +A\pd_zv) =0
\eeq
if and only if   $v$   satisfies \rea{bva} and
\ga\label{deq1}
\pdoz v+b+A\pd_zv=0. 
\end{gather}
\item Let $u$ be
a continuous map from $\ov D$ into an open subset $\Om$ of $\cc^n$ with
$v\in\cL C^1(D)$. Let $A\in\cL C^1(\Om)$ be an $n\times n$ matrix.
Then    $v$ satisfies
\eq{ieq2}
v +\T (b +A(v){\ov \pd_zv}) =0
\eeq
if and only if    $v$    satisfies \rea{bva} and
\ga\label{deq2}
\pdoz v+b+A(u)\ov{\pd_zv}=0. 
\end{gather}
\eppp
Here equations \rea{bva}-\rea{deq2} are on   $D$.
\end{lemma}


We prove a version of \rt{regid-} with parameter.
\pr{regid}    Let $0<\all<1$ and let $k\geq j\geq0$ be   integers.
 Let $\Om$ be a bounded domain in $\cc$ with $\pd\Om\in\cL C^{k+1+\all}$.    Let
 $a,b\in\B^{k+\all,j}(\ov\Om,P)$ be (scalar) functions
satisfying  $
|a |_{\all,0}<\e_\all$.
Then
\eq{dfdg++}
 v(\cdot,t)+\T_{\Om} b(\cdot,t)+\T_{\Om}(a(\cdot,t)\pd_zv(\cdot,t))=0
\eeq
has a unique solution $v(z,t)$ with $v\in\B^{k+\all+1,j} (\ov\Om,P)$.
Consequently,
$$\I+Ta\pd_z \colon \B ^{k+1+\all,j}(\ov\Om,P)\to
\B ^{k+1+\all,j}(\ov\Om,P)$$
 has a bounded inverse.
\epr
\begin{proof}
By \rl{isa}, there exists a solution $v\in\B^{1+\all,j} $
to \re{dfdg++}. The proposition is verified  for $k=0$. If $k\geq1$, the assertion that $v\in\B^{k+1+\all,j} $
follows from \rp{combv+} and \rl{teqs1}.
The last assertion   in the proposition follows from $T_\Om(\B^{k+\all,j}(\ov\Om,P))=\B^{k+1+\all,j} (\ov\Om,P)$
and the open mapping theorem.
\end{proof}

\begin{rem}\label{hatcb+}
Let $0<\all<1$ and $k,j$ be nonnegative integers.
Let $\Om$ be a bounded domain in $\cc$ with $\pd\Om\in\cL C^{k+1+\all}$.
It would be interesting to know if
\eq{}
\nonumber
\I+Ta\pd_z \colon  \hat{\cL C}^{k+1+\all,j}(\ov\Om,P)\to
\hat{\cL C}^{k+1+\all,j}(\ov\Om,P)
\eeq
has  a bounded inverse,
assuming $a\in\hat{\cL C}^{k+\all,j}(\ov\Om,P)$ and
$\|a\|_{\all,0}$ is small.  \rl{isa}
is for the case when $a(\cdot,t)$ has compact support.
\end{rem}
%

We now use the  proof of \rp{combv+}   to study a problem in   different directions.
However, unlike the previous case, the next one fails
   in higher dimension.
\pr{probc}  Let $0<\all<1$ and $k \geq0$ be an  integer.
Let $\gaa$ be an embedded curve in $\cc$ of class
$\cL C^{k+1+\all}$. Let $\Om_1,\Om_2$ be disjoint open subsets of $\cc$ such that
  both  $\pd\Om_1,\pd\Om_2$ contain
$\gaa$ as relatively open subsets. Assume that
 $a_i \in\B^{k+\all}(\Om_i\cup\gaa)$
 satisfies $
|a_i|_{0,0}<1$ on $\Om_i\cup\gaa$. Let $E$ be an embedded $\cL C^1$ curve in  $\D$ such that
$\D\setminus E$ is open in $\cc$ and
 has exactly two connected components $\om_1,\om_2$. Assume that $u$
is a continuous map from $\D$ into $\Om_1\cup \gaa\cup \Om_2$ such that $u\colon\om_i\to\Om_i$
are $J$-holomorphic with respect to $\pdoz+a_i\pd_z$. Then $E$ is a    curve of class
$\cL C_{loc}^{k+1+\all}$.
\epr
\begin{proof} The proof is a slight modification of  the proof of \rp{combv+}.
The problem is local. Fix $z_0\in E$ and let $p=u(z_0)$. We may assume that near $p$, $\gaa$ is contained in
the real axis and $\Om_1,\Om_2$ are contained in the lower and upper half planes. Applying
a local
  change of coordinates $\varphi_i$ which is of class $\cL C^{k+1+\all}$
   on $\Om_i\cup\gaa$ and fixes $\gaa$ pointwise,
we may assume that
$a_j(x,y)=O(|y|^{k+\all})$. Let $a $   be $a_i$ on $\Om_i\cup\gaa$.
Then $X=\pdoz+a\pd_z$ is of class $\cL C^{k+\all}$ on $\Om_1\cup \gaa\cup\Om_2$.
Near  $p\in \gaa$, we apply a diffeomorphism $\phi$ of class $\cL C^{k+1+\all}$ which transforms
$X$ into  $\pdoz$. Let $g=\phi\circ\phi_i\circ u$ on $\om_i\cup E$, which is holomorphic away from $E$.
Since $g$ is continuous and $E$ is an embedded
$\cL C^1$ curve, then  $g$ is holomorphic at $z_0$. It is easy to verify that
  $g$ is biholomorphic near $z_0$. Consequently, $E$ is of class $\cL C^{k+1+\all}$
near $z_0$.
\end{proof}

\begin{exmp}
Let $E$ be an embedded
 $\cL C^1$ curve    connecting $i, -i$ and dividing  $\D$ into two components $\om_1,\om_2$.
Let $\la$ be a $\cL C^\infty$  function on $\D$ which is positive on $\om_1$
and negative on $\om_2$. The existence of such a function is trivial, by taking
it   vanishing to infinity order along $E$. We use the standard
complex structure on $\D\times\cc=
\{\IM w<\la(z)\}\cup\{\IM w=\la(z)\}\cup
\{\IM w>\la(z)\}$. Let $u(z)=(z,0)$. Then $u\colon\om_i\to\Om_i$
are holomorphic, $\gaa=\{\IM w=\la(z)\}$ is $\cL C^\infty$, $u(\D)=\D\times 0$,
but $E$ needs not be
$\cL C^\infty$.
\end{exmp}

\smallskip
The main purpose of next result is to provide another proof of \rp{combv+}. The proof
does not yield a sharp result. To get the sharp result, we have to return to the
argument in \rp{combv+}.
We will deal with non-tangential boundary values. We will restrict to the non-parameter case.

\pr{combv++}   Let $0<\all<1$ and $k\geq j\geq0$ be   integers.
Let $\gaa$ be an embedded curve in $\cc$ of class
$\cL C^{k+1+\all}$. Let $\Om_1,\Om_2$ be disjoint open subsets of $\cc$ such that
  both  $\pd\Om_1,\pd\Om_2$ contain
$\gaa$ as relatively open subsets. Assume that
 $a_i,b_i\in\B^{k+\all}(\Om_i\cup\gaa)$
 satisfy $
|a_i|_{0,0}<1$ on $\Om_i\cup\gaa$.
Suppose  that $f|_{\Om_i}$ are continuous and admit the same non-tangential boundary value
function $f\in L^p(\gaa)$ with $p>1$.
Let $f$   satisfy
\eq{pdfab}
\pdoz f+a_i \pd_zf=b_i\ \text{on $\Om_i$},\qquad i=1,2.
\eeq
 Then   $f\in\B_{loc}^{k+1+\all}
(\Om_i\cup\gaa)$.
\epr
\begin{proof}
We may assume that $\Om_1,\Om_2$ are two disjoint bounded simply connected
domains whose
boundaries are of class $\cL C^{k+1+\all}$.
Apply a $\cL C^{k+1+\all}$ diffeomorphism $\psi_i$
of $\ov\Om_i$ onto $\ov\Om_i'$
which transforms $\pdoz+a_i\pd_z$ into $\pdoz$. Such   $ \psi_i$ exists in view of  \rp{isop} by
 extending $a_i$ to a neighborhood of $\gaa$ via Whitney's extension theorem and by shrinking
 $\Om$ at $p\in\gaa$.
 By a theorem of Kellogg,
there exists a Riemann mapping $ \phi_i\in\cL C^{k+1+\all}(\ov{ \Om_i'})$ which
sends $\Om_i'$ onto the upper half-plane. We may assume that $\gaa\neq\pd\Om_j$ and $\gaa$
is mapped into a compact subset by $\phi_j\circ\psi_j$.

Without loss of generality, we may assume that $\gaa=(-1,1)$. We choose subdomain $\om_j$ of $\Om_j$
as follows:
  $\pd\om_j$ contains $
[-r_0,r_0]$;   $f$ has
non-tangential limits at $r_0, -r_0\in\gaa$;
$\phi_j\circ\psi_i$  sends $\ov\om_j$ onto
  $Q=[r',r'']\times[0,1]$.
Now, let $\phi$ be a Riemann mapping for $Q$. Note that $\phi$ is smooth on $\ov Q$ and
$\phi'=0$, $\phi''\neq0$ at vertices  of $Q$.
Let $\var_j=\phi\circ\phi_j\circ\psi_j$.
Thus,  $(\var_j^{-1})^*L^p(\pd\om_j)\subset L^p(\pd\D)$.
(For our local results,   we avoid the use of $\var_j^*(L^p(\pd\D))\subset L^q(\pd\om_j)$ for
 $q<p/2$.)

Let $\Ht$ be the    conjugate
operator on $\pd\D$. Namely, for a real function $f\in L^p(\pd\D)$ with $1<p<\infty$,
there is a holomorphic function
$h$ on $\D$ with $\IM h(0)=0$
 whose non-tangential boundary value is $f+i\Ht f$ with $\Ht f$ real-valued
 (Theorem 3.1,   p.~57; Lemma 1.1, p.~103 in~\ci{Gaeion}).
 The $\mathcal A_if=(\Ht (f\circ\varphi_i^{-1}))\circ\varphi_i$ is called
the  conjugate operator on $\pd\om_j$ for $\pdoz+a_i\pd_z$.

By a lemma of M.~Riesz  (\cite{Gaeion}, p.~113), $\|\Ht v\|_{L^p(\pd\D)}\leq C_p\|v\|_{L^p(\pd\D)}$
for $1<p<\infty$. Thus,
 $\Ht^2 f=-f+c_f$ for   $f\in L^p(\pd\D)$ with $c_f$ being a constant.
By Privalov's theorem,
$\Ht (L^p(\pd\D)\cap\cL C^{k+\all}(E))\subset  \cL C_{loc}^{k+\all}(E)$ for an arc $E$ in $\D$.

From now on, we assume that $1<p<\infty$.
By our choice of $\om_j$, $f$ is bounded on $\pd\om_j\setminus \gaa$.
Thus $f|_{\om_j}$ has   non-tangential limit functions     in $L^p(\pd\om_j)$.
 Recall that for $u\in L^p(\pd\D)$,
$$ 
\Ht u(z)= -\f{1}{\pi}\, p.v.\int_{\pd \D}u(\zeta)\, d\log|\zeta-z|,\quad z\in\pd\D.
$$ 
Thus for $u_i\in L^p(\pd\om_j)$,
\al\label{aba}
\mathcal A_i g(z)&=- \f{1}{\pi}\, p.v.\int_{\pd\om_i }u_i(\zeta)\, d\log|\var_i (\zeta)-\var_i (z)|\\
&\df - \f{1}{\pi}\lim_{\e\to0}\int_{\pd\om_i\cap\{|\var_i(\zeta)-\var_i(z)|>\e\}}u_i(\zeta)\, d\log|\var_i (\zeta)-\var_i (z)|.
\nonumber
\end{align}
Let $f_i  \in\cL C^{k+1+\all}(\ov\Om_i )$ be a solution to  the inhomogeneous equation   \re{pdfab}. We get
\eq{eqh}
f=f_i+g_i,\quad \text{on $\ov\Om_i $},\quad \pdoz g_i+a_i\pd_zg_i=0.
\eeq
Our assumption  implies that $g_i|_{\om_j}$ has non-tangential limits in $L^p(\pd\om_j)$.
By \re{aba}-\re{eqh}, we have
$$
g_i=u_i+\sqrt{-1}\mathcal A_iu_i+ic_i.
$$
Let $\chi_\gaa$ be the characteristic function of $\gaa$. Obviously,
 $E_i=\mathcal A_i((1-\chi_\gaa)u_i)\in\cL C_{loc}^{k+1+\all}(\gaa)$.
On $\gamma$, we have $f_1+g_1=f_2+g_2$ and hence
$$
u_1+\RE f_1=u_2+\RE f_2, \quad \mathcal A_1u_1+ \IM f_1+c_1=\mathcal A_2 u_2+\IM f_2+c_2.
$$
  By the first identity, we obtain
$\mathcal A_2 (\chi_\gaa(u_2-   u_1))\in\cL C_{loc}^{k+1+\all}(\gaa)$. The second   shows
$$
\mathcal A_2(\chi_\gaa u_1)   -\mathcal A_1(\chi_\gaa u_1)\in\cL C_{loc}^{k+1+\all}(\gamma).
$$
We assume that $\Om_2$ and $\gaa$ have the same orientation.   On $\gaa$, we rewrite \re{aba} as
\aln
\mathcal A_1 (\chi_\gaa u_1)(z)&=\f{1}{\pi}\, p.v.\int_{\gaa}
u_1 (\zeta)  \, d\log|\var_1(\zeta)-\var_1(z)|,
\\
\mathcal A_2 (\chi_\gaa u_1)(z)& = -\f{1}{\pi}\, p.v.\int_{ \gaa}u_1(\zeta)\, d\log|\var_2( \zeta) - \var_2(z)|.
\end{align*}
Here
the change of sign arises from the opposite orientation of $\gaa$ in $\pd\om_1$.

Assume now that $k\geq1$. We may assume that $\gaa=(0,1)$. Therefore, on $\gaa$
\aln
(\mathcal A_1  -\mathcal A_2)(\chi_\gaa u_1)(x)&= 2\mathcal A_1(\chi_\gaa u_1)(x)
+E_3(x)+C(x),
\end{align*}
where
\aln E_3(x)
&=\f{1}{\pi}\int_{0}^{1} u_1(t) \, d
\log\f{
|\var_2(t)-\var_2(x)|}{|\var_1(t)-\var_1(x)|}
+C(x),\\
C(x)&= \lim_{\e\to0}\biggl\{
\int_{I_2(x,e)}-\int_{I_1(x,e)}\biggr\}
\f{u_1 (t)}{\pi}\,
  d\log|\var_1( t) - \var_1(x)|.
\end{align*}
Here $I_j(x,\e)=(0,1)\setminus (x-\e_j,x+\e_j')$  with $e_j, e_j'>0$, and
   $$|\var_j(x-(-1)^j\e_j)-\var_j(x)|=
|\var_j(x+(-1)^j\e_j')-\var_j(x)|=
c_j(x)^{-1}\e.$$
Also, $c_j(x)>0$ satisfies
$|\varphi_j(t)-\varphi_j(x)|=c_j(x)|t-x+O(|t-x|^{1+\all})|$. Note that
 $e_j=e+O(\e^{1+\all})$
and $e_j'=e+O(e^{1+\all})$.
This shows that $(I_2(x,\e)\setminus I_1(x,\e))\cup(I_1(x,\e)\setminus I_2(x,\e))$ is
contained in $[e-Ce^{1+\all}, e+Ce^{1+\all}]\cup [-e-Ce^{1+\all}, -e+Ce^{1+\all}]$.
  Therefore,   we can verify that
\aln
C(x)&=\f{1}{2\pi}\lim_{\e\to0}\biggl\{
\int_{I_2(x,\e)}-\int_{I_1(x,\e)\}}\biggr\}\,
\f{u_1 (t)}{x-t}\, dt=0
\end{align*}
when $x$ is a Lebesgue point of
$u_1\in  L^1_{loc}(\gaa)$.  Hence, $  C=0$ a.e.~on $\gaa$.
Since $E_3\in\cL C_{loc}^{k-1+\all}$, we obtain
$
2\mathcal A_1 (\chi_\gaa u_1)  \in\cL C_{loc}^{k-1+\all}(\gaa).
$
Hence, $\Ht ((\chi_\gaa u_1)\circ\var_1^{-1})\in L^p(\pd\D)\cap\cL C_{loc}^{k-1+\all}(\tilde\gaa)$
for $\tilde\gaa=\var_1(\gaa)$.
Since $(\chi_\gaa u_1)\circ\var_1\in L^p(\pd\D)$,  then $\Ht^2f=-f+c_f$
implies $(\chi_\gaa u_1)\circ\var_1^{-1}\in
L^p(\pd\D)\cap\cL
C_{loc}^{k-1+\all}(\tilde\gaa)$. Therefore,
$$g_1\circ\var_1^{-1}|_{\pd\D}=
u_1\circ\var_1^{-1}+i\Ht(u_1\circ\var_1^{-1})+ic_1\in L^p(\pd\D)\cap
\cL C_{loc}^{k-1+\all}(\tilde\gaa).
$$
 And  $g_1\circ\var_1^{-1}\in
\cL C_{loc}^{k-1+\all}(\D\cup\tilde\gaa).
$
 Hence $g_1\in \cL C^{k-1+\all}_{loc}(\om_1\cup\gaa)$.
Note that when $k\geq2$, we can apply integration
 by parts and achieve  $E_3\in\cL C^{k+\all}_{loc}$, and  hence $g_1\in\cL
 C_{loc}^{k+\all}
 (\Om_1\cup\gaa)$.

\medskip

The above argument does not yield the sharp result.
We now turn to the proof by using  previous methods.   We need to use
a Fatou lemma.
 We   assume that $\gaa=(-1,1)$.
As used above, the $g_j$, given by \re{eqh}  and holomorphic
in $\pdoz+a_j\pd_z $,
has a non-tangential limit function in $L_{loc}^p(\gaa)$. We choose $\om_j,\var_j$ as before. So
$g_j\circ\var_j^{-1}$ is holomorphic on $\D$ with non-tangential limit function     in $L^p(\pd\D)$.
 By Fatou's lemma,
$g_j\circ\var_j^{-1}(re^{i\theta})-g_j\circ\var_j^{-1}(e^{i\theta})$
and hence
  $f \circ\var_j^{-1}(re^{i\theta})-f \circ\var_j^{-1}(e^{i\theta})$ tend to zero
in $L^p(\pd\D)$
 as $r\to1^-$. Assume that $\gaa=(-1,1)$. Let $r_0$ be given in the definition of $\om_j$. Fix $0<r_1<r_0$.
Write $[-r_1,r_1]=\{\var_j^{-1}(e^{i\theta}); \theta_j\leq
\theta\leq\theta_j'\}$.
Let $\gaa_{j,t}=\{\var_j^{-1}(te^{i\theta})\colon \theta_j\leq
\theta\leq\theta_j'\}$. Then $f|_{\gaa_{j,t}}$ tends to $f|_{\gaa_{j,1}}$ in $L^p$ norm as $t\to1^-$. More precisely,
$\int_{\theta_j}^{\theta_j'}|f(\gaa_{2,t}(\theta))-f(\gaa_{2,1}(\theta))|^p\, d\theta$ tends to $0$ as $t\to1^-$.
Note that $\gaa_{1,1},\gaa_{2,1}$ are the same set $[-r_1,r_1]$ with opposite orientations.

As in the proof of \rp{combv+}, taking smaller $r_1$ if necessary and using
three changes of coordinates
$\phi_1,\phi_2,\psi$ and a solution of an inhomogeneous equation,
we arrive at the case that $a_j, b_j$ are zero.
Now, $f$ is holomorphic
away from $\tilde\gaa=\psi(\gaa)$,  $f|_{\gaa_{j,t}}$ tends to
$f|_{\tilde\gaa_{j,1}}$ in $L^p$ norm as $t\to1^-$,
and $\tilde\gaa_{1,1}, \tilde\gaa_{2,1}$ are the same curve with opposite orientations.
 Applying  the Cauchy formula
to  cancel boundary integrals in $\tilde\gaa_{1,1},\tilde\gaa_{2,1}$,
we find the extension of
$f$   to a neighborhood of $p\in\tilde\gaa_{1,1}$ via a Cauchy transform on  a small circle
  centered at $p$. Returning to the original
coordinates, we  obtain the desired conclusion for
the original $f$.
\end{proof}

\medskip

As mentioned in the introduction, our main result fails for harmonic functions.
Let $\Om$ be a bounded domain in $\cc$ with $\pd \Om\in \cL C^{\infty}$.
 Suppose that $f$ is continuous on $\pd\Om$ and $dt$ is the arc-length element on $\pd\Om$.
Then
$
W_f(z)=
\f{1}{\pi } \int_{\pd\Om} f(t)\log|\gaa(t)-z|\, dt
$
  is harmonic on $\cc\setminus\pd\Om$ and continuous on $\cc$.
 However,
\aln
\pd_{n(s)}W_f&= f(s)+\f{1}{\pi}\int_{\pd\Om} f(t)\pd_{s}\arg (\gaa(s)-\gamma(t))\, dt, \\
\quad\pd_{-n(s)}W_f&=  f(s)-\f{1}{\pi}\int_{\pd\Om} f(t)\pd_{s}\arg (\gaa(s)-\gamma(t))\, dt.
\end{align*}
Here $n(t)$ is the unit outer normal vector of $\pd\Om$.  In particular, if $f$ is not
smooth, then $W_f$ cannot be  smooth simultaneously on $\ov\Om$ and $\cc\setminus \Om$.

We  remark that if $W_f\in\cL C^1(\cc)$, then $f$ and
 $W_f$
must be  zero.


\newcommand{\NWsith}{\bibitem{NWsith}
A. Nijenhuis and W.B. Woolf,
{\it Some integration problems in almost-complex and complex manifolds},
Ann. of Math. (2) {\bf 77}(1963) 424--489.}

\newcommand{\IRzefo}{\bibitem{IRzefo}
S. Ivashkovich and J.-P. Rosay, {\em
 Schwarz-type lemmas for solutions of
 $\overline\partial$-inequalities and complete
 hyperbolicity of almost complex manifolds},
 Ann. Inst. Fourier (Grenoble)  {\bf 54}(2004),  no. 7, 2387--2435.}

\newcommand{\Befise}{\bibitem{Befise}L. Bers,
{\it Riemann Surfaces} (mimeographed lecture notes), New York
University, (1957-1958). }

\newcommand{\Vesitw}{\bibitem{Vesitw}
I.N. Vekua,  {\it Generalized analytic functions},
 Pergamon Press, London-Paris-Frankfurt; Addison-Wesley
Publishing Co., Inc., Reading, Mass. 1962. }

\newcommand{\Fonifi}{\bibitem{Fonifi} G.B. Folland, {\it
Introduction to partial differential equations},
Second edition. Princeton University Press, Princeton, NJ, 1995.
}

\newcommand{\HJeini}{\bibitem{HJeini}
N. Hanges and H. Jacobowitz, {\em  A remark on almost complex structures with boundary},
Amer. J. Math.  {\bf 111}(1989),  no.~1, 53--64.}

\newcommand{\GTzeon}{\bibitem{GTzeon}
D. Gilbarg and N.S. Trudinger, {\em Elliptic partial differential equations of second order}.
Reprint of the 1998 edition. Classics in Mathematics. Springer-Verlag, Berlin, 2001.}

\newcommand{\Runion}{\bibitem{Runion}
W. Rudin, {\it Functional analysis}, Second edition.
International Series in Pure and Applied Mathematics. McGraw-Hill, Inc., New York, 1991.
}

\newcommand{\NNfise}{\bibitem{NNfise}A. Newlander and L. Nirenberg, {\em Complex
analytic coordinates in almost complex manifolds},  Ann. of Math. (2) {\bf
65}(1957), 391--404. }

\newcommand{\Weeini}{\bibitem{Weeini}
S.M. Webster, {\em A new proof of the Newlander-Nirenberg theorem}, Math. Z. {\bf
201}(1989), no. 3, 303--316. }

 \newcommand{\Caeiei}{\bibitem{Caeiei}D.W. Catlin,
 {\em  A Newlander-Nirenberg theorem for manifolds with boundary},
Michigan Math. J.  {\bf 35}(1988),  no. 2, 233--240. }

\newcommand{\Sesifo}{\bibitem{Sesifo}
R.T. Seeley, {\em Extension of $C\sp{\infty }$ functions defined
in a half space}, Proc. Amer. Math. Soc.  {\bf 15}(1964),
625--626. }

\newcommand{\Hieiei}{\bibitem{Hieiei}
C.D. Hill, {\em
What is the notion of a complex manifold with a smooth boundary?}
Algebraic analysis, Vol.~I, 185–-201, Academic Press, Boston, MA, 1988.
}

\newcommand{\Honize}{\bibitem{Honize}
L. H\"ormander, {\em The analysis of linear partial differential
operators.
 I. Distribution theory and Fourier analysis.}
 Springer-Verlag, Berlin, 1990.}

\newcommand{\Gaeion}{\bibitem{Gaeion} J.B. Garnett, {\it Bounded analytic functions},
 Pure and Applied Mathematics, 96. Academic Press, Inc., New York-London, 1981.
 }

\end{document}